\documentclass[a4paper,12pts]{amsart}
\usepackage{graphicx} 
\usepackage{lineno}

\usepackage{amsthm}
\usepackage{amssymb}
\usepackage{amsfonts}
\usepackage{mathtools}
\usepackage{mathrsfs}
\usepackage{graphicx}
\usepackage{color}
\usepackage{enumitem}
\usepackage{setspace}
\usepackage[english]{babel}
\usepackage{csquotes}
\usepackage{etoolbox}
\usepackage{thmtools}
\usepackage{yfonts}
\usepackage{mfirstuc}
\usepackage{stix}
\usepackage{fancyvrb}
\usepackage{textgreek} 

\usepackage{amsmath}
\usepackage{amsthm, amssymb, amsfonts}
\usepackage{graphicx, color}
\usepackage{caption, subcaption} 
\usepackage{multicol} 
\usepackage{stmaryrd}
\usepackage{tikz-cd} 
\usepackage{lipsum}  
\usepackage[capitalise,noabbrev, nameinlink]{cleveref} 
\Crefname{part}{Part}{Parts}
\Crefname{step}{Step}{Steps}
\Crefname{prop}{Proposition}{Propositions}
\Crefname{prob}{Problem}{Problems}
\newcommand{\forces}[1]{\Vdash\text{\say{#1}}}
\newcommand{\Lim}{\text{Lim}}
\newcommand{\Fin}{\text{Fin}}
\newcommand{\type}{\text{type}}

\usepackage{enumitem}
\renewcommand{\leq}{\leqslant}
\renewcommand{\geq}{\geqslant}

\renewcommand{\epsilon}{\varepsilon}
\newcommand{\M}{\mathcal{M}}
\newcommand{\E}{\mathbb{E}}
\newcommand{\X}{\mathbb{X}}
\newcommand{\A}{\mathcal{A}}

\newcommand{\stickT}{%
\setbox255=\hbox{\raise1ex\hbox{$\hspace{0.2pt}\,\bullet\,$}}
\mathord{\rlap{\hbox to\wd255{\hss\hbox{$|$}\hss}}
\box255}
}
\newcommand{\stickS}{%
\setbox255=\hbox{\raise0.6ex\hbox{$\scriptstyle\bullet$}}
\mathord{\rlap{\hbox to\wd255{\hss\hbox{$\scriptstyle|$}\hss}}
\box255}
}

\numberwithin{equation}{section}
\theoremstyle{plain}
\newtheorem{theorem}[equation]{Theorem}
\newtheorem{lemma}[equation]{Lemma}
\newtheorem{proposition}[equation]{Proposition}
\newtheorem{corollary}[equation]{Corollary}

\theoremstyle{definition}
\newtheorem{definition}[equation]{Definition}

\newtheorem{problem}[equation]{Problem}

\theoremstyle{remark}

\usepackage{mathtools}
\usepackage{amsfonts}
\usepackage{amssymb}
\usepackage{amsthm}
\usepackage{amsmath}
\usepackage{dirtytalk}
\usepackage{mathrsfs}  
\usepackage{lineno}

 \usepackage{multicol}  
 \usepackage{bbm}
 \usepackage{tabularx}
\newenvironment{claimproof}[1][\proofname]
{\begin{proof}[#1]}
{\end{proof}}

\title{There may be an $n$-entangled set but no $n+1$-entangled sets}
\author{Jorge Cruz}
\begin{document}
\maketitle
\begin{abstract}In this paper we answer Question 3.5 of \cite{carroy2025questionsentangledlinearorders} in the positive by showing  that for every $2\leq n\in \mathbb{N}$, the statement \say{there is an $n$-entangled set, but there are no $n+1$-entangled sets} is consistent. We also prove some theorems which improve our understanding of entangled sets in relation to construction schemes: (1) The axiom FCA$^\Delta$ introduced in \cite{finitizationclubch} implies the existence of $n$-entangled sets which are not $n+1$-entangled. (2) $\mathfrak{m}_\mathcal{F}>\omega_1$ implies the non-existence of entangled sets. Thus, $2$-capturing schemes alone are not sufficient to build these kinds of linear orders. (3) The existence of a 2-$\Delta$-capturing scheme is consistent with MA.  
\end{abstract}
\footnotetext{Keywords: entangled, construction schemes, morasses, set theory, CH, Baumgartner axiom. Mathematics Subject Classification (2020): 03E65, 03E75, 03E50, 03E05, 03E35.}\footnotetext{The research on tis paper is partially supported by grants from NSERC(455916) and CNRS(UMR7586). The Author declares no conflict of interest}.The topic of this work has its roots in  \cite{baumgartneraxiom}. There, Baumgartner showed the consistency, via PFA, of the axiom which is now known as Baumgartner's axiom BA$(\omega_1)$ stating that all $\aleph_1$-dense set of reals are isomorphic. Abraham and Shelah later complemented Baumgartner's work by showing in \cite{MAdoesnotImplyBA} that the negation of $BA(\omega_1)$ is consistent with MA. In order to prove this result they introduced the concept of an entangled set of reals. Given $1\leq n\in \omega$ we say that an uncountable set $\E\subseteq \mathbb{R}$ is \emph{$n$-entangled} provided that for every $\tau:n\longrightarrow 2$ (referred to as a \emph{type}) and each uncountable  family $\mathcal{A}\subseteq [\E]^n$ whose elements are pairwise disjoint, there are distinct $a,b\in \mathcal{A}$ so that $a(i)<b(i)$\footnote{Here, $a(i)$ and $b(i)$ stand for the ith elements of $a$ and $b$ with respect to their increasing order.} if and only if $\tau(i)=0$. $\E$ is said to be entangled whenever it is $n$-entangled for each $n$. The study of $BA(\omega_1)$ and entangled sets continued in the very influential \cite{abrahamrubinshelah}. Among other things,  Abraham, Rubin and Shelah proved there that $BA(\omega_1)$ is consistent with $\mathfrak{c}:=2^\omega>\omega_2$, and that the statement \say{There is an $n$-entangled set for every $1\leq n\in \omega$, but there are no entangled sets} is consistent as well, even in the presence of MA. This last result can be used to produce models of MA in which the Magidor-Malitz language is not countably compact. Recently, Guzm\'an and Todorc\v{e}vi\'c improved in \cite{guzmanentangled} the result of \cite{MAdoesnotImplyBA} by proving the consistency of the statement: \say{There is a $2$-entangled set}+MA+PID. In particular, this means that the existence of a $2$-entangled set of reals is not equivalent to an equality of the form \say{cardinal invariant=$\omega_1$} in the presence of PID, at least for those cardinals whose value is decided as $\mathfrak{c}$ by MA. This result contrasts with the ones in \cite{CombinatorialDichotomiesandCardinalInvariants} or \cite{GapsandTowers} where it is shown that under PID, the existence of certain uncountable objects with strong combinatorial properties (a sixth Tukey type and a non-Hausdorff tower respectively) is equivalent to an equality between certain cardinal invariants and $\omega_1.$\\
The study of Bamgartner's axiom, as well as some of its modifications and/or generalizations, is still active. Beyond the papers mentioned above, it is worth to highlight \cite{Cdhmartinsaxiom}, \cite{WhyYcc}, \cite{schemescruz}, \cite{BAomega2mooretodorcevic}, \cite{homeomorphismsofmanifolds} and \cite{PartitionProblems}. One of the many things proved by Todor\v{c}evic  in this last famous book is that BA$(\omega_1)$ implies that $\mathfrak{b}>\omega_1.$ In \cite{switzer2025weakbaumgartneraxiomsuniversal}, Switzer provides a more complete discussion of the state of the art regarding this subject, including both  famous open problems as well as posing some new interesting ones.
As far as entangled sets is concern, there has also been further development of the theory, for example \cite{AsperoMotaEntangled}. In \cite{RemarksonChainConditionsinProducts}, Todor\v{c}evi\'c  shows that the existence of an entangled set implies that being $ccc$ is not a productive property and that $cof(\mathfrak{c})=\omega_1$ implies the existence of an entangled set. He also extends the definition of $n$-entangled sets to arbitrary linear orders and sizes and observe that $2$-entangled orders are  $ccc$ and that $3$-entangled orders are separable (thus, they can be thought as subsets of $\mathbb{R}$). Due to this facts, it follows that any $2$-entangled linear order which is not separable, is necesarilly a Suslin line. The consistency of the existence of such a line was first obtained by Krueger in \cite{entangledsuslinkrueger} via forcing, where he also extendes the definition of entangledness to arbitrary \emph{partial} orders. More recently in \cite{carroy2025questionsentangledlinearorders}, Carroy, Levin and Notaro adressed some  questions regarding this topic.  In that same paper, they ask:
\begin{problem}[Question 3.5, \cite{carroy2025questionsentangledlinearorders}]Let $n\geq 2$. Can there be an $n$-entangled set but no $(n+1)$-entangled sets?
    
\end{problem}
This question is very natural, and in some sense, a positive answer to it (which we provide in this paper),  helps to justify the study of $n$-entangled sets iniciated in \cite{MAdoesnotImplyBA}  for the complete spectrum of all natural numbers instead of just working with $2$-entangled sets and entangled sets.   Our main theorem is the following:
\begin{theorem}\label{maintheorem}
    Let $n\geq 2$. The statement \say{There is an $n$-entangled set, but there are no $n+1$-entangled sets}+PID+MA is consistent.
\end{theorem}
In order to prove this theorem, we will reformulate the results of \cite{guzmanentangled}, where they implicitly prove the consistency  of  the  forcing axiom associated to all proper forcings preserving the $n$-entangledness of a fixed $n$-entangled set $\E$ . Having this in mind, our main job will be to take an arbitrary uncountable set $\X$ and produce a proper forcing which destroys the $n+1$-entangledness of $\X$ while preserving the $n$-entangledness of $\E.$ 

In \cite{carroy2025questionsentangledlinearorders}, it is  shown that under CH there is, for each $n>1$, an $n$-entangled set which is not $n+1$-entangled, and that there is an entangled set $\mathbb{E}$  which is homeomorphic to a subset of $\mathbb{R}$ which is not even $2$-entangled. In this paper we show that the axiom FCA$^\Delta$ as defined in \cite{finitizationclubch} suffices to prove both of these theorems. We remark that $FCA^\Delta$ follows from CH and holds in the Cohen model. It is also worth pointing out that we provide a single proof for both theorems and ours is singificantly shorter.
\begin{theorem}[Under $FCA^\Delta$]\label{schemenotn+1theorem} There is a family $\langle \mathbb{E}_n\rangle_{n\in \omega}$ of pairwise homeomorphic uncountable subsets of $\mathbb{R}$ so that $\mathbb{E}_0$ is entangled, $\mathbb{E}_1$ is not $2$-entangled and $\mathbb{E}_n$ is an $n$-entangled but not $n+1$-entangled set for any $n\geq 2.$
\end{theorem}
The axiom FCA$^\Delta$ is associated to the objects known as \emph{construction schemes}\footnote{All undefined notions related to schemes will be defined later in Section \ref{preliminariesschemesection}.} introduced for the first time by Todor\v{c}evi\'c in \cite{schemenonseparablestructures} and further studied in papers such as  \cite{schemescruz}, \cite{irredundantsetsoperator}, \cite{banachspacescheme}, \cite{lopezschemethesis} and \cite{treesgapsscheme}. We finish this paper by proving two results regarding the limitations of construction schemes. Namely, we show that the existence of a $2$-$\Delta$-capturing schemes is consistent with MA and that the existence of a $2$-capturing scheme is consistent with the non-existence of entangled sets. The main importance of these results lies in the fact that the existence of $2$-$\Delta$-capturing or $2$-capturing scheme already imply the the failure of BA$(\omega_1).$ Formally, we show:
\begin{theorem}[Under $\mathfrak{m}^\Delta_\mathcal{F}>\omega_1$]\label{mfdeltatheorem} Assume $\mathbb{P}$ is a $ccc$-forcing. Then $\mathbb{P}$ $2$-$\Delta$-preserves $\mathcal{F}$.
\end{theorem}
\begin{theorem}[Under $\mathfrak{m}_\mathcal{F}>\omega_1$]\label{mFtheorem}There are no entangled sets. 
\end{theorem}
In view of the results presented here, the author is left with the following problem.
\begin{problem} Fix $2\leq n\in \omega$. Is it consistent that there is an $n$-entangled (resp. entangled) and every $2$-entangled set is $n$-entangled (resp. entangled)?
\end{problem}

The structure of the paper is as follows: In \emph{Section} \ref{notationsection} we introduce some general notation. In \emph{Section} \ref{preliminariesentangledsection} we present some preliminaries regarding entangled sets. In \emph{Section} \ref{maintheorem} we prove Theorem \ref{maintheorem}. In \emph{Section} \ref{preliminariesschemesection} we provide a brief introduction on construction schemes.  In \emph{Section} \ref{schemenotn+1section} we prove Theorem \ref{schemenotn+1theorem}. Finally, in \emph{Section} \ref{mfsection} we prove theorems \ref{mFtheorem} and \ref{mfdeltatheorem}.

\section{Notation}\label{notationsection}
The notation and terminology used here is mostly standard and it follows \cite{schemescruz}.  Given a set $X$ and a (possibly finite) cardinal $\kappa$, $[X]^\kappa$ denotes the family of all subsets of $X$ of cardinality $\kappa$. The sets $[X]^{<\kappa}$ and $[X]^{\leq \kappa}$ have the expected meanings. The family of all non-empty finite sets of $X$ is denoted by $\Fin(X)$. That is, $\Fin(X)=[X]^{<\omega}\backslash\{\emptyset\}$. $\mathscr{P}(X)$ denotes the power set of $X$. By $\text{Lim}$  we mean the  set of limit ordinals strictly smaller than $\omega_1$. Given a partial order $(\X,<)$ and $X,Y\subseteq \X$, we write $X<Y$ whenever $x<y$ for each $x\in X$ and $y\in Y$. In the same way, if $\mathcal{X},\mathcal{Y}\subseteq \mathscr{P}(X)$, then $\mathcal{X}<\mathcal{Y}$ if and only if $X<Y$ for each $X\in \mathcal{X}$ and $Y\in \mathcal{Y}.$ Given $A,B\subseteq \mathbb{X}$, $A$ is said to be an \emph{initial segment} of $B$, denoted as $A\sqsubseteq B$, if for any $a\in A$ and $b\in B$, if $b<a$, then $b\in A$. In this same situation we also say that $B$ \emph{end-extends} $A.$

For a well-ordered set $X$ (typically finite), we denote by $ot(X)$ its order type. We identify $X$ with the unique strictly increasing function $h:ot(X)\longrightarrow X$. In this way, $X(\alpha)=h(\alpha)$ denotes the $\alpha$ element of $X$ with respect to its increasing enumeration. Analogously, $X[A]=\{X(\alpha)\,:\,\alpha\in A\}$ for $A\subseteq ot(X)$.  A family $\mathcal{D}$ is called a \textit{$\Delta$-system} with root $R$ if $|\mathcal{D}|\geq 2$ and $X\cap Y= R$ whenever $X,Y\in \mathcal{D}$ are different.  If moreover, $R<X\backslash R$ for any and $X\backslash R<Y\backslash R$ or viseversa for any two $X,Y\in \mathcal{D}$, we call $\mathcal{D}$ a root-tail-tail $\Delta$-sytem. 

Recall that if $f,g:\omega\longrightarrow\omega$ are distinct,  we can define the number $\Delta(f,g)$ as the minimum $k$ for which $f(k)\not=g(k).$ By convenience, we put $\Delta(f,g):=\omega$ whenever $f=g$.

Given sets $U_{0},\dots, U_{n-1}$, we define $U_0\otimes U_1:=\{u\cup v\,:\,u\in U_0\text{ and }v\in U_1\}$ and more generaly $\otimes_{i<n} U_i:=\{\bigcup_{i<n} u(i)\,:\,u\in \prod\limits_{i<n}U_i\}$. for  $n\in \omega$, $i\in 2$ a function $\tau:n\longrightarrow 2$, we define $\tau^i:n\longrightarrow 2$ as $\tau^i(j)=|i-t(j)|$. Given $n,m\in \omega$, $\sigma_0:n\longrightarrow 2$ and $\sigma_1:m\longrightarrow 2$, we define the concatenation of $\sigma_0$ and $\sigma_1$ as the function $\sigma_0\frown\sigma_1:n+m\longrightarrow 2$ given by: $$\sigma_0\frown \sigma_1(i)=\begin{cases}\sigma(i)&\text{ if }i<n\\
\sigma(i-n)&\text{ if }i\geq n
    \end{cases}
$$

Whenever we say that $M$ is an elementary submodel of $H(\kappa)$, we mean that $(M,\in )$ is an elementary submodel of $(H(\kappa),\in )$, and whenever we say $\kappa$ is a large enough cardinal, we mean that $H(\kappa)$ contains every object of interest and satisfies enough ZFC axioms to carry the arguments discussed in the respective proof.

Given $n\geq 1$, we inherit $[\mathbb{R}]^n$ with the \emph{Vietoris topology} whose base consist of all sets of the form $\otimes_{i<n} U_i$ where $U_i$ is open in $\mathbb{R}$ for each $i<n$, and $U_i\cap U_j=\emptyset$ for any two distinct $i,j<n$.  It is well-known and easy to prove that $[\mathbb{R}]^n$ is homeomorphic to $\{(a_0,\dots,a_{n-1})\in \mathbb{R}^n\,:\,a_0<\dots<a_{n-1}\}$. In particular, it is second countable.

Given a forcing notion $\mathbb{P}$, we denote by $\mathfrak{m}(\mathbb{P})$ its Martin's number, i.e. , the minimal cardinal $\kappa$ for which there are $\kappa$-many dense sets over $\mathbb{P}$ which are not intersected by any filter. The cardinal $\mathfrak{m}$ denotes the minimum of all $\mathfrak{m}(\mathbb{P})$ where $\mathbb{P}$ ranges over all $ccc$-forcing notions. For us, MA denotes the statement $\mathfrak{m}=\mathfrak{c}>\omega_1$. 
\section{Some preliminaries on entangled sets}\label{preliminariesentangledsection}
This small section is mainly devoted to fix some notation regarding entangled sets as well as review some easy but important results regarding these objects. For a more  detailed treatment we refer the reader to \cite{guzmanentangled}.\\
Let $1\leq n \in \omega$ and $a,b\in [\mathbb{R}]^n$ be disjoint sets. We define $\type(a,b):n\longrightarrow 2$ given by: $$\type(a,b)(i)=\begin{cases}0&\text{ if }a(i)<b(i)\\
1&\text{ if }a(i)>b(i)
\end{cases}$$
In the context of entangled sets, we refer to a function $\tau:n\longrightarrow 2$ as a \emph{(configuration) type}. We say that $a$ and $b$ have \emph{type $\tau$} provided that $\type(a,b)=\tau$. Note that $\type(a,b)$, in case it is well-defined, is always distinct from $\type(b,a)$. In fact,  $\type(a,b)=\type(b,a)^1.$ In order to take into consideration this fact, we define for an arbitrary $\tau$ as above, $\Vert \tau \Vert=\{\tau^0,\tau_1\}$. In this way, $\Vert \type(a,b)\Vert=\Vert \tau\Vert$ if and only either $\type(a,b)=\tau$ or $\type(a,b)=\tau^1.$ The following lemma is well-known and easy to prove.
\begin{lemma}Let $2\leq n$ and $\mathbb{E}\in[ \mathbb{R}]^{\omega_1}$. Then $\mathbb{E}$ is $n$-entangled if and only if for each $\tau:n\longrightarrow 2$ and every uncountable family $\{a_\alpha\,:\,\alpha\in \omega_1\}\subseteq [\mathbb{R}]^n$  whose elements are pairwise disjoint, there are $\alpha<\beta\in \omega_1$ so that $\type (a_\alpha,a_\beta)=\tau$.   
\end{lemma}

The definition of the $\type$ function can be further extended. Assume that $U,V\subseteq \mathscr{P}([\mathbb{R}]^n)$ are non-empty and such that $\bigcup U\cap \bigcup V=\emptyset$. If $\tau:n\longrightarrow 2$ is such that $\type(a,b)=\tau$ for any two $a\in U$ and $b\in V$, we define $\type(U,V)$ as $\tau$.  Whenever we assert that $\type(U,V)=\tau$, we will implicitly assume that $U$ and $V$ satisfy the requirements mentioned above. Furthermore, if $a\in  [\mathbb{R}]^n$ and $V\subseteq \mathscr{P}([\mathbb{R}]^n)$,then $\type(\{a\},V)$ is simply written as $\type(a,V).$

The following proposition is well-known. It implicitly says that \say{being disjoint and having type $\tau$} defines an open  graph over $[\mathbb{R}]^n.$ We will frequently make use of this property without any explicit mention.
\begin{proposition} Assume that $a,b\in [\mathbb{R}]^n$ are disjoint. Then there are open sets $U,V\subseteq  [\mathbb{R}]^n$ so that $\bigcup U\cap \bigcup V=\emptyset$, $a\in U$, $b\in V$ and $$\type(a,b)=\type(U,V).$$ 
\end{proposition}

The next two lemmas already appear in some form in \cite{abrahamrubinshelah} and \cite{MAdoesnotImplyBA}. Both of them are heavily used in the \say{duplication} type of arguments presented in \cite{abrahamrubinshelah}. The first one is just a consequence that both $\min$ and $\max$ are definable. The second one was independently proven by Todor\v{c}evi\'c.

\begin{lemma}[Duplication lemma for total orders]\label{duplicationorder} Assume that $\X\subseteq \mathbb{R}$ is an uncountable set and $M$ is a countable elementary submodel of a large enough $H(\kappa)$ so that $\X\in M$. Then for any $x\in \X\backslash M$, there are $y,z\in \X\cap M$ so that $y<x<z.$
    
\end{lemma}
\begin{lemma}[Duplication lemma for entangled sets]\label{duplicationentangled} Let $2\leq n\in \omega$. Assume that  $\E\subseteq \mathbb{R}$ an $n$-entangled set,  $\mathcal{A}\subseteq [\E]^n$ is uncountable and $M$ is a countable elementary submodel of a large enough $H(\kappa)$ so that $\A\in M$. Then for each type $\tau:n\longrightarrow 2$ and every $a\in \mathcal{A}$ disjoint from $M$, there is $b\in \mathcal{A}\cap M$ (necessarily disjoint from $a$) so that $\type(a,b)=\tau.$
\end{lemma}
The main result of of \cite{guzmanentangled} was the consistency of MA+PID+\say{There is a $2$-entangled set}. In order to achieve this result, Guzm\'an  and  Todor\v{c}evic proved an preservation theorem for Neeman iterations (as defined by Neeman in \cite{nemmanmodelstwotypes} and $2$-entangled sets. It is easy to see that their argument  actually serves to prove a much more general preservation theorem for uncountable graphs without uncountable independent sets. In particular, such preservation theorem can be written in terms of $n$-entangled sets as follows.
\begin{theorem}[Theorem 94, \cite{guzmanentangled}]Let $2\leq n \in \omega$,  $\theta$ be an inaccesible cardinal, $J, \mathcal{S},\mathcal{T},\mathbb{P}^{\mathcal{S},\mathcal{T}}_\in$ and $\mathbb{P}$ be such that:
\begin{itemize}
    \item $\mathcal{S}\subseteq\{M\in [H(\theta)]^\omega\,:\,M\text{ is elementary submodel of }H(\theta)\}$ is stationary in $[H(\theta)]^\omega$.
    \item $\mathcal{T}=\{H(\lambda)\,:\,H(\lambda)\text{ is elementary submodel of }H(\theta)\text{ and }cof(\lambda)>\omega\}.$
    \item $\mathbb{P}^{\mathcal{S},\mathcal{T}}_\epsilon$ is the two type side condition $\epsilon$-collapse associated to $\mathcal{S}$ and $\mathcal{T}$ (see Definition 58 of \cite{guzmanentangled}).
    \item $J:\theta\longrightarrow H(\theta)$.
    \item $\mathbb{P}=\mathbb{P}(J)$ is the Neeman iteration associated with $J$ and $\mathbb{P}^{\mathcal{S},\mathcal{T}}_\epsilon$ (see Definition 72 of \cite{guzmanentangled}).
   
\end{itemize}
Let $\mathbb{E}$ be an $n$-entangled set. If for every $X\in \mathcal{T}$ either $X$ is trivial or $$\mathbb{P}\cap X\forces{$J((X)\text{ preserves the $n$-entangledness of }\E$ },\footnote{preserving the $n$-entangledness of $\mathbb{E}$ and forcing $\mathbb{E}$ to be $n$-entangled is the same thing.}$$ then $\mathbb{P}\forces{$\E$\text{ is }$n$\text{-entangled}}.$
\end{theorem}
As a direct consequence of this, we get a parametrized version of PFA for a fixed $n$-entangled set. Let us first explicitly describe this new axiom.\\\\
\noindent
\textbf{PFA$_{n-ent}(\mathbb{E})$:} $\mathbb{E}$ 
is an $n$-entangled set and for any proper forcing $\mathbb{P}$, if $\mathbb{P}$ preserves the $n$-entangledness of $\mathbb{E}$, then $\mathfrak{m}(\mathbb{P})$.\\\\
\noindent
The same argument proved in Theorem 97 of \cite{guzmanentangled} yields the consistency of PFA$_{n-ent}(\mathbb{E})$. This as well as a general proper forcing axiom associated to families a graphs without independent sets will appear in an upcomming paper which is joint work with Guzm\'an and Todor\v{c}evi\'c. Instances of these parametrized forcing axioms have been studied  independently by Erlebach (see \cite{emilythesis}).

\begin{theorem}[Esentially Theorem 97,\cite{guzmanentangled}]Let $2\leq n\in \omega$ and $\mathbb{E}$ be an $n$-entangled set. The is an extension of the universe via a proper forcing in which $\mathfrak{c}=\omega_2$ and PFA$_{n-ent}(\mathbb{E})$ holds.
\end{theorem}
In Theorem 37 of \cite{guzmanentangled} it was implicitly proved that $\mathbb{P}_{2-ent}(\mathbb{E})$ implies MA. It turns out that this result holds for arbitrary $n$'s. Even though the proof exposed in \cite{guzmanentangled} can be modified to yield the following theorem, we present a different proof here since ours is just a few lines long.
\begin{theorem}$\mathbb{P}_{n-ent}(\mathbb{E})$ implies $\mathfrak{m}>\omega_1$.  
\begin{proof}Let $\mathbb{Q}$ be a $ccc$ forcing.  Let $\mathbb{P}_0$ be $\sigma$-centered forcing collapsing $\mathfrak{c}$ to $\omega_1$. Note that $\mathbb{P}_0$ preserves the $n$-entangledness of $\mathbb{E}$ and the fact that $\mathbb{Q}$ is $ccc$.  In this way, $\mathbb{P}_0\times \mathbb{Q}$ is a proper forcing. If it turns out that $\mathbb{P}_0\times \mathbb{Q}$ preserves the $n$-entangledness of $\mathbb{E}$, we are done. This is because $\mathfrak{m}(\mathbb{P}_0\times \mathbb{Q})\leq \mathfrak{m}(\mathbb{Q})$. On the other hand, if $\mathbb{P}_0\times \mathbb{Q}$ does no preserve the $n$-entangledness of $\mathbb{E}$, it is necessarily true that $$\mathbb{P}_0\forces{ $\mathbb{Q}$\text{ does not preserve the $n$-entangledness of $\mathbb{E}.$ } } $$

This is due to the fact that $n$-entangled sets are preserved by $\sigma$-centered forcings. Now, we evoke  the claim 1 of Lemma 8.5 in \cite{abrahamrubinshelah}: \\\\
\noindent
\underline{Claim 1}(\cite{abrahamrubinshelah}): Assume CH holds and let $\mathbb{Q}$ be a ccc-forcing which does not preserve the $n$-entangledness of $\mathbb{E}$. Then there is a $ccc$-forcing $\mathbb{P}_1$ which preserves the $n$-entangledness of $\mathbb{E}$ and which adds an uncountable antichain to $\mathbb{Q}.$\\\\
Since $\mathbb{P}_0$ forces CH, we can fix a $\mathbb{P}_0$-name $\dot{\mathbb{P}}_1$ for a $ccc$-forcing satisfying the conclusion of the claim above. In order to finish, just note that $\mathbb{P}_0\ast \dot{\mathbb{P}}_1$ is a proper forcing which preserves the $n$-entangledness of $\mathbb{E}$ and adds an uncountable chain to $\mathbb{Q}$. As PFA$_{n-ent}(\E)$ holds and being an uncountable antichain is a property which can be coded by $\omega_1$-many dense sets, we conclude that $\mathbb{Q}$ is not $ccc$. However, this is a contradiction.  
\end{proof}
\end{theorem}
Theorem 47 of \cite{guzmanentangled} can also be generalized to $n$-entangled sets. We do not prove such generalization since the arguments are identical to the original ones.
\begin{theorem}[Essentially Theorem 47, \cite{guzmanentangled}]PFA$_{n-ent}(\E)$ implies PID.
\end{theorem}

\section{Destruction of $n+1$-entangledness}
The main purpose of this section is to prove Theorem \ref{maintheorem}. In order to do this, it suffices to show that PFA$_{n-ent}(\E)$ implies the non-existence of $n+1$-entangled sets. Let us fix $2\leq n\in \omega$, an $n$-entangled set $\mathbb{E}$,  an uncountable $\mathbb{X}\subseteq \mathbb{R}$ and $\kappa=(2^\mathfrak{c})^+$. We define $Sm(\kappa)$ as the set of all countable elementary submodels of $H(\kappa)$ having $\mathbb{E}, \mathcal{B}$ and $\mathbb{X}$ as elements. We say that $\mathcal{M}\in [Sm(\kappa)]^{<\omega}$ is a \emph{chain} whenever it is totally ordered with respect to $\in$. In this situation, we write the elements of $\mathcal{M}$ increasingly as $\mathcal{M}(0),\dots,\mathcal{M}(|\M|-1)$. By convention, $\mathcal{M}(|\M|)$ denotes the universe.

Our goal is to define a proper forcing $\mathbb{P}(\X,\E)$ which preserves the $n$-entangledness of $\E$ and which force $\X$ to not be $n+1$-entangled. The key of our proof lies in fixing the two unique constant types $\overline{0}$ and $\overline{1}$ over $[\X]^n$ as the ones we want to avoid,  and construct via finite approximations an uncountable set $\mathcal{X}\subseteq[\X]^n$ such that the type of their elements is not constant. In order for this approach to work, it will be necessary that the elements of each $x\in \mathcal{X}$ are separated by certain chains of elementary submodels.

In the definition below, the $x$'s will latter be interpreted as possible candidates for elements of the set $\mathcal{X}\subseteq [\X]^n$ which will be forced by $\mathbb{P}(\X,\E)$ to avoid the constant types.

\begin{definition}[Good pairs] We say that $(\mathcal{M},x)\in [Sm(\kappa)]^{n+1}\times [\X]^{n+1}$ is a \emph{good pair} provided that $x(i)\in \M(i+1)\backslash \M(i)$ for each $i\leq n$. Given any other good pair $(\mathcal{N},y)$, we write $(\M,x)\in _\star(\mathcal{N},y)$ if $(\M,x)\in \mathcal{N}(0)$. The set of all good pairs is denoted as $Gp(\kappa)$.
\end{definition}
Note that $\in_\star$ is partial order over $Gp(\kappa)$.  Whenever we talk about chains, initial segments, etc., inside $Gp(\kappa)$ we will be referring to this order. We are know ready to define the required forcing notion.
\begin{definition}[The forcing $\mathbb{P}(\X,\E)$] We define $\mathbb{P}(\X,\E)$ to be the set of all finite $\in_\star$ chains $p\subseteq G_p(\kappa)$ with the extra property that if $(\mathcal{M},x),(\mathcal{N},y)\in p$ are distinct, then $\type(x,y)$ is not constant. We order $\mathbb{P}\X$ by reverse inclusion. Given any $p\in \mathbb {P}(\X,\E)$, we let $n_p=|p|$ and write the elements of $p$ increasingly as $(\mathcal{M}^p_0,x^p_0)\in_\star\dots\in_\star(\M_{n_p-1}^p,x_{n_{p-1}})$.
\end{definition}
We omit the proof of the following easy lemma.
\begin{lemma}\label{lemmasumodel}Assume that $p\in \mathbb{P}(\X,\E)$ and $M$ is an elementary submodel of $H(\kappa)$ with $p\in M$. Then there is a good pair $(\mathcal{M},x)$ so that $p\cup\{(\mathcal{M},x)\}\in \mathbb{P}(\X,E)$ and $\mathcal{M}(0)=M$.
\end{lemma}
The proof of the theorem below can be thought as an easier version of the proof of the preservation of the $n$-entangledness of $\E$ which will be presented below in Theorem \ref{mainpreservationentangledness}. For that reason, we leave it to the reader.
    
\begin{theorem}Let $\theta$ be a large enough cardinal so that $\mathbb{P}(\X,\E),\kappa\in H(\theta)$. and let $\hat{M}$ be an elementary submodel of $H(\theta)$ with $\mathbb{P}(\X,\E),\kappa\in \hat{M}$. If $p\in \mathbb{P}(\X,\E)$ is a condition so that $\hat{M}\cap H(\kappa)=\mathcal{M}_i(0)$ for some $i<n_p$, then $p$ is $(\hat{M},\mathbb{P}(\X,\E))$-generic.
\end{theorem}
As a direct corollary of the two results above we get:
\begin{corollary} $\mathbb{P}(\X,\E)$ is a proper forcing.  Furthermore, if $G\subseteq \mathbb{P}(\X,\E)$ is a generic filter, then $\mathcal{X}_G:=\{x^p_i\,:p\in G\text{ and }i<n+1\}$ is forced to be an uncountable subset of $[\X]^{n+1}$ whose elements are pairwise disjoint and so that $\type(a,b)$ is not constant for any two distinct $a,b\in \mathcal{X}_G$. In particular, $\X$ is not $n+1$-entangled in $V[G].$
\end{corollary}

The next definition together with the two lemmas below it are the key for proving the preservation of the $n$-entangledness of $\E$ via the forcing $\mathbb{P}(\X,\E)$.
\begin{definition}Let $\mathcal{M}\in [Sm(\kappa)]^{n+1}$ and $a\in [\E]^{\leq n}$. Given $i\leq n$ we define $$a^{i\M }:=a\cap (\M(i+1)\backslash \M(i)),$$
We say that:
\begin{itemize}
    \item $a$ is of the \emph{first kind} over $\mathcal{M}$ if $|a^{i\mathcal{M}}|=1$ for each $i<n$.\footnote{Note that in this case $|a|=n$.}
    \item $a$ is of the \emph{second kind} over $\mathcal{M}$ if it is not of the first kind.
    
\end{itemize}
\end{definition}
One main obstacle for proving the preservation of the $n$-entangledness of $\E$ is that at some point we will want to duplicate points in the second coordinates of good pairs, as well as $\leq n$-size subsets of $\mathbb{E}$ via lemmas \ref{duplicationorder} and \ref{duplicationentangled} respectively. Unfortunately, we can not apply these two lemmas at the \say{same time}. The main motivation for defining sets of the first and second kind is to address this problem. Naively speaking if a set is of the first kind, we will prioritize duplications via Lemma \ref{duplicationorder}. On the other hand, if we are dealing with a set of the second kind, we will prioritize duplications via Lemma \ref{duplicationentangled}. As it can bee seen in the next two lemmas, the main job of the chains of submodels forming the first coordinate of good pairs is to leave enough room of error for our arguments to work.

\begin{lemma}\label{lemmakind1}Let $(\M,x)\in Gp(\kappa)$, $a\in [\E]^{ n}$ be of the first kind over $\mathcal{M}$, and $\tau:n\longrightarrow 2$ be an arbitrary type. If $\mathcal{D}\in \M(0)\cap \mathscr{P}([X]^{n+1}\times[\mathbb{E}]^n)$ is such that $(x,a)\in \mathcal{D}$, then there is $(y,b)\in \mathcal{D}\cap \M(0)$ so that:\begin{enumerate}
    \item  $\type(x,y)$ is not constant,
    \item  $\Vert\type(a,b)\Vert=\Vert \tau\Vert$.
\end{enumerate}
\begin{proof}  For the sake of simplicity we will assume that $a^{i\M}<a^{j\M}$ whenever $i<j\leq n$.  Note that in this situation $a^{i\M}$ coincides with the set whose only element is the ith element of $a$ with respect to its increasing order, namely $a(i)$. Furthermore, we will assume that $\tau(n-1)=0$. We proceed to build by reverse recursion two families of open subsets of $\mathbb{R}$ $\langle U_i\rangle_{i\leq n}$ and $\langle V_i\rangle_{i<n}$ such that the following properties are satisfied:
\begin{enumerate}[label=(\Roman*)]
    \item $U_0<\dots<U_{n-1}<U_n,$
    \item $V_0<\dots<V_{n-1}$,
        \item $\type(x[\{n-1,n\}], U_{n-1}\otimes U_n)=\langle0,1\rangle$. In other words, $x(n-1)>U_{n-1}$ and $x(n)<U_n$.
    \item $\Vert\type(a, \,\otimes_{i<n} V_i)\Vert=\Vert \tau\Vert$.
\end{enumerate}
\noindent
\emph{Step $n$ (Base step):} Consider the set $\mathcal{D}_n=\{y\in X\,:\,x(n-1)<y\text{ and }(x[n]\cup\{y\},a)\in \mathcal{D}\}$. Note that $x(n)\in \mathcal{D}_n$, and since $\mathcal{D}\in \M(0)\subseteq \mathcal{M}(n)$, it also holds that $\mathcal{D}_n\in\mathcal{M}(n)$.  According to the definition of good pair,  $x(n)\not\in \mathcal{M}(n)$. In this way, we can find a point $y\in \mathcal{D}_n\cap  \mathcal{M}(n)$ so that $x<y$. We define $U_n$ to be an open set inside $\mathcal{M}(0)$ such that $y\in U_n$ and $x(n)<\overline{U_n}.$
\\\\
\noindent
\emph{Step $n-1$(Second step):} Consider the set $\mathcal{D}_{n-1}$ of all $y\in \X$ so that $x(n-2)<y<\overline U_n$ and for which there are $z\in \mathbb{E}$ and $u\in [U_n]^1$ satisfying the following properties:
\begin{itemize}
    \item $a(n-2)<z,$
    \item$(x[n-1]\cup\{y\}\cup u,a[n-1]\cup\{z\})\in \mathcal{D}$.
    \end{itemize}
Note that $x(n-1)\in \mathcal{D}_{n-1}$ and $\mathcal{D}_{n-1}\in \M(n-1) $. Since $x(n-1)\not\in \M(n-1)$, we can find $y\in \mathcal{D}_{n-1}\cap \mathcal{M}(n-1)$ so that $y<x(n-1)$. Let $z\in \mathbb{E}\cap \mathcal{M}(n-1)$ witness the fact that $y\in \mathcal{D}_{n-1}$. Observe that $z\not=a(n-1)$ because $a(n-1)$. We now define $U_{n-1}$ and $V_{n-1}$ to be open subsets of $\mathbb{R}$ inside $\mathcal{M}(0)$ so that:
\begin{itemize}
    \item $y\in U_{n-1}$ and $x(n-2)<\overline{U_{n-1}}<x(n-1)$,
    \item $z\in V_{n-1}$ and $a(n-2)<\overline{V_{n-1}}$,
    \item $a(n-1)<\overline{V_{n-1}}$ or $\overline{V_{n-1}}<a(n-1)$. In this way, $\type(a(n-1),V_{n-1})$ is well-defined.\\
\end{itemize}
We have now completed the first two steps of the recursion. Before we continue, let us define $\epsilon\in 2$ as $0$ if $a(n-1)<\overline{V_{n-1}}$, and as $1$ otherwise.\\\\
\noindent
\emph{Recursive step}: Assume that $k+1\leq {n-1}$ and we have defined $V_{k+1},\dots, V_{n-1}, U_{k+1},\dots,U_{n}\in \mathcal{B}$ in such way that the following recursive hypotheses are fulfilled: \begin{enumerate}
    \item $a(k)<\overline{ V_{k+1}}<\dots<\overline{V_{n-1}}$,
    \item $x(k)<\overline{U_{k+1}}<\dots<\overline{U_{n-1}}<\overline{U_n}$,
    \item $\type(x[\{n-1,n\}], U_{n-1}\otimes U_n)=\langle0,1\rangle$,
    \item $\type\big(a[n\backslash (k+1)], \: \bigotimes\limits_{k+1\leq i<n}V_i\,\big)(j)=\tau^\epsilon((k+1)+j)$ for each $j<n-(k+1),$
    \item There are $u\in \bigotimes\limits_{k< i\leq n} U_i$ and  $v\in \bigotimes\limits_{k<i<n} V_i$ so that $(x[k]\cup u,a[k]\cup v)\in \mathcal{D}.$
\end{enumerate}
Let us define $\mathcal{D}_k$ as the set of all $z\in \E$ so that\footnote{In here, both $a(k-1)$ and $x(k-1)$ denote $-\infty$ in the case where $k=0.$ In this same situation, both $x[k-1]$ and $a[k-1]$ denote the empty set.} $a(k-1)<z<\overline{V_{k+1}}$ and for which there are  $u\in \bigotimes\limits_{k< i\leq n} U_i$,  $v\in \bigotimes\limits_{k<i<n} V_i$ and $y\in \mathbb{X}$ satisfying the following properties:
\begin{itemize}
\item $x(k-1)<y<\overline{U_{k+1}},$
\item $(x[k-1]\cup \{y\}\cup u,a[k-1]\cup\{z\}\cup v)\in \mathcal{D}.$
\end{itemize}
Observe that $a(k)\in \mathcal{D}_k$ and $\mathcal{D}_k\in \M(k)$. Since $\{a(k)\}=a^{i\mathcal{M}}$, we also have that $a(k)\not\in \mathcal{M}(k)$. In this way, we can find a point $z\in \mathcal{D}_k\cap \M(k)$ so that $a(k)<z$ if $\tau^\epsilon(k)=0$ and $z<a(k)$ if $\tau^\epsilon(k)=1.$ Let $y\in \X\cap \M(k)$ be a witness of the fact that $z\in \mathcal{D}$. We now define $U_k$ and $V_k$ to be two open sets inside $\mathcal{M}(0)$ so that:\begin{itemize}
    \item $y\in U_k$ and $x(k-1)<\overline{U_k}<\overline{U_{k+1}},$
    \item $z\in  V_k$ and $z(k-1)<\overline{V_k}<\overline{V_{k+1}},$
    \item $a(k)<\overline{V_k}$ if $a(k)<z$ and $a(k)>\overline{V_k}$ if $z>a(k).$
\end{itemize}
By defining $U_k$ and $V_k$ in this way, it is easy to see the recursive hypotheses of $k$ are fulfilled. This finishes the recursion.\\\\
Note that $U_0,\dots ,U_n,V_0,\dots,V_{n-1}$ have been constructed in such way that the properties (I), (II), (III) and (IV) are satisfied. Furthermore, by construction we have that $$\big((\bigotimes_{i<n+1} U_i)\times( \bigotimes_{i<n} V_i)\big)\cap \mathcal{D}\not=\emptyset.$$
As both $\mathcal{D}$ and $\mathcal{B}$ are elements of $\M(0)$, there is an pair $(y,b)$ in the intersection above which also belongs to $\M(0)$.  It is clear that such $(y,b)$ finishes the job.
\end{proof}

\end{lemma}

\begin{lemma}\label{lemmakind2}Let $(\M,x)\in Gp(\kappa)$, $a\in [\E]^{\leq n}$ be of the second kind over $\mathcal{M}$, and $\sigma:|a|\longrightarrow 2$ be an arbitrary type. If $a\cap \mathcal{M}(0)=\emptyset$ and $\mathcal{D}\in \M(0)\cap [\X]^n\times[E]^{|a|}$ is such that $(x,a)\in \mathcal{D}$, then there is $(y,d)\in \mathcal{D}\cap M(0)$ so that:\begin{enumerate}
    \item  $\type(x,y)$ is not constant,
    \item  $\type(a,b)= \sigma$.
\end{enumerate}
\begin{proof}Again, we will assume for simplicity that $a^{i\M}<a^{j\M}$ whenever $i<j\leq n.$  We also let $P:=\{i\leq n\,:\,a^{i\M}\not=\emptyset\}$ and for each $i\leq n+1$ we let $m_i\leq|a|$ be such that $a[m_i]=a\cap \mathcal{M}(i)$. Note that if $i\in P$, then $m_i=\min(m<|a|\,:\,a(m)\in a^{i\M})$. We proceed to define a family  $\langle U_i\rangle_{i\leq n}$ of open subsets of $\mathbb{R}$ and family $\langle V_i\rangle_{i\in P}$ with $V_i$ is an open subset of $[\mathbb{R}]^{|a^{i\M}|}$ for each $i\in P$, and so that the following properties hold:
\begin{enumerate}[label=(\Roman*)]
    \item $U_0<\dots<U_n$,
    \item $V_i<V_j$ whenever $i<j\in P$,
    \item $\type(x,\otimes_{i\leq n} U_i)$ is not constant,
    \item $\type (a, \otimes_{i\in P}V_i)=\tau$.
\end{enumerate} 
\emph{Step $n$ (Base step):}  Here, we shall consider the two following cases:\\\\
\noindent 
\emph{Case 1:} If  $n\not\in P$ we simply find $U_n$ the same way as we did in the base step of Lemma \ref{lemmakind1}. That is, we find  an open  $U_n\in \mathcal{M}(0)$ so that $x(n)<\overline{U_n}$ and there is $y\in \X\cap U_n$ with $(x[n-1]\cup\{y\},a)\in \mathcal{D}$.\\\\
\noindent
\emph{Case 2:} f $n\in P$, we make use of Lemma \ref{duplicationentangled} to find a point $z\in [\mathbb{E}]^{|a^{n\M}|}\cap \mathcal{M}(n)$ so that:\begin{itemize}
    \item $a[k_n]<z$,
    \item $\type(a^{n\M},z)(j)=\sigma(m_n+j)$ for each $j<|a^{n\M}|$,
    \item there is $x(n-1)<y\in \X $ distinct from $x(n)$ such that $(x[n-1]\cup\{y\},a[m_n]\cup z)\in \mathcal{D}.$
\end{itemize}    
    We then proceed to take an open $V_n\subseteq [\mathbb{R}]^{|a^{n\M}|}$ and an open $U_n\subseteq \mathbb{R}$, both inside $\mathcal{M}(0)$, such that $z\in V_n$, $x\in U_n$, $a[k_n]<\overline{V_n}$, $x(n-1)<\overline{U_n}$, $a^{n\M}\cap \bigcup V_n=\emptyset$, $\type(a^{n\M},V_n)=\type(a^{n\M},z)$,  $x(n)<U_n$ if $x(n)<y$, and $x(n)>U_n$ otherwise.\\\\
    \noindent
\emph{Recursive step}: Assume that $k+1\leq {n}$ and we have defined $\langle U_i\rangle_{k<i\leq n}$ and $\langle V_i\rangle_{k<i\in P}$ all of whose elements are in $\mathcal{M}(0)$ and in such way that the following recursive hypotheses are fulfilled: \begin{enumerate}
    \item $x(k)<\overline{U_{k+1}}<\dots<\overline{U_{n-1}}<\overline{U_n}$,
    \item $a^{i\M}<V_j$ for any $i<j\in P$,
    \item for any $k<j\leq n$, either $x(j)<\overline{U_j}$ or $\overline{U_j}<x(j)$,
    \item for each $k<j\in P$ and every $i<|a^{j\M}|$, $\type(a^{j\M},V_j)(i)=\sigma(m_j+i)$.
    \item there are $u\in \bigotimes\limits_{k< i\leq n} U_i$ and  $v\in \bigotimes\limits_{k<i\in P} V_i$ so that $(x[k]\cup u,a[m_{k+1}]\cup v)\in \mathcal{D}.$
\end{enumerate}
We again proceed by cases:\\\\
\emph{Case 1:} If $k\not\in P$, note that $a[m_{k+1}]=a[m_k]$. Here, we simply use Lemma  \ref{duplicationorder} to find an open $U_k\in \mathcal{M}(0)$ satisfying the following properties:
\begin{itemize}
\item $x(k-1)<\overline{U_k}<\overline{U_{k+1}}$,
    \item $x(k)<\overline{U_k}$ if $x(k+1)>\overline{U_{k+1}}$ and $x(k)>\overline{U_k}$ otherwise,
    \item there are $y\in U_k$,  $u\in \bigotimes\limits_{k< i\leq n} U_i$ and  $v\in \bigotimes\limits_{k<i\in P} V_i$ so that $(x[k-1]\cup\{y\}\cup u,a[m_k]\cup v)\in \mathcal{D}.$
\end{itemize}
\noindent
\emph{Case 2:} If $k\not\in P$ we use Lemma \ref{duplicationentangled} to find two open sets $V_k\subseteq [\mathbb{R}]^{|a^{k\M}|}$ and $U_k\subseteq \mathbb{R}$ inside $\mathcal{M}(0)$ with the following properties:
\begin{itemize}
    \item $x(k-1)<\overline{U_k}<\overline{U_k+1}$,
    \item $a^{k-1\M}<\overline{V_k}< V_{i_{k+1}}$ where $i_{k+1}=\min(i\in P\,:\,i>k)$,
    \item either $x(k)<\overline{U_k}$ or $\overline{U_k}<x(k)$,
    \item  $a^{k\M}\cap \bigcup V_k=\emptyset$ and $\type(a^{k\M},V_k)(i)=\sigma(m_j+i)$ for each $i<|a^{k\M}|.$
    \item There are $y\in U_k$, $z\in V_k$, $u\in \bigotimes\limits_{k< i\leq n} U_i$ and  $v\in \bigotimes\limits_{k<i\in P} V_i$ so that $(x[k-1]\cup\{y\}\cup u,a[m_k]\cup\{z\}\cup v)\in \mathcal{D}.$
\end{itemize}
It is easy to see that the recursive hypotheses hold for $k$ independently of the case. This finishes the construction.\\
It is straightforward from points $(1)$, $(2)$, and $(4)$ of the recursive hypotheses that the sequences $\langle U_i\rangle_{i\leq n}$ and $\langle V_i\rangle_{i\in P}$ satisfy the properties (I), (II) and (IV). The key observation is that since $a$ is of the second kind over $\mathcal{M}$, then there is at least one $k<n$ such that $k\in P$. For this $k$, we followed the Case 1 of the recursive step when constructing $U_k$.  According to the second item in such case, we have that $\type(x[\{k,k+1\}],U_k\otimes U_{k+1})$ is not constant. From this, it easily follows that the property (III) is satisfied. In order to finish, just note that the item (5) from the recursive hypotheses assures that there is $y\in \otimes_{i<n} U_i$ and $b\in \otimes_{i\in P} V_i$ so that $(y,b)\in \mathcal{D}$. The existence of this pair completes the proof.
\end{proof}
\end{lemma}
A particular instance of the lemma above happens when $a=\sigma=\emptyset$. In such case, we can rewrite it in the following manner.
\begin{lemma}\label{lemmakind2empty} Let $(\mathcal{M},x)\in G_p(\kappa)$ and $\mathcal{D}\in \mathcal{M}(0)\cap[X]^n$ be such that $x\in\mathcal{D}$. Then there is $y\in \mathcal{D}\cap M(0)$ such that $\type(x,y)$ is not constant.
    
\end{lemma}
The proof of Theorem \ref{maintheorem} will conclude once we prove the following theorem.

\begin{theorem}\label{mainpreservationentangledness}$\mathbb{P}(\X,\E)$ preserves the $n$-entangledness of $\E.$
\begin{proof} Let $\dot{\A}$ be a name for an uncountable subset of $[\E]^n$ whose elements are pairwise disjoint and let $q_0\in \mathbb{P}(\X,\E)$. Fix $\theta$ a large  enough cardinal so that $\mathbb{P}(\X,\E),\kappa\in H(\theta)$. Now, let $\hat{M}$ be an elementary submodel of $H(\theta)$ having $\mathbb{P}(\X,\E)$, $\dot{\A}$ and $q_0$ as elements. By applying Lemma \ref{lemmasumodel}, we can find  $x\in [\X]^{n+1}$, and $p\leq q$  so that $(M,x)\in p$. Since $\dot{\A}$ is forced to be uncountable and its elements are forced to be pairwise disjoint, we may assume without any loss of generality that there is $a\in [\E]^n$ disjoint from $M$ for which $p\forces{ $a\in \dot{A} $ }$.\\
Given $i<n_p$, let us write $\mathcal{M}^p_i$ and $x^p_i$ simply as $\mathcal{M}_i$ and $x_i$ respectively. Also, fix $r:=M\cap p$.  Note that $r$ is a condition end-extended by $p$. Furthermore, $\mathcal{M}_{n_r}(0)=M.$  Our plan is to find a condition $q\in \mathbb{P}(\X,\E)\cap M$ and $b\in [\mathbb{E}]^n$ so that:
\begin{itemize}
    \item $q$ end-extends $r$ and its compatible with $p$,
    \item $q\forces{ $a\in \dot{A}$ }$,
    \item $\Vert\type(a,b)\Vert=\Vert \tau\Vert.$
\end{itemize}
In order to do this, we first fix $n_{rp}:=n_p-n_r$ and for each $i<$ we let\footnote{Here,  $\mathcal{M}_{n_p}(0)$ is understood as the whole universe.} $$a_i:=(\mathcal{M}_{n_r+i+1}(0)\backslash \mathcal{M}_{n_r+i}(0))\cap a.$$
For the sake of simplicity we will assume that $a_i<a_j$ whenever $i<j<n_{rp}.$ Now, let  $P=\{i<n_{rp}\:\,a_i\not=\emptyset\}$ and fix two sequences $\langle O_i\rangle_{i<{n_{rp}}}$ and $\langle W_i\rangle_{i\in P}$ inside $M$ such that the following properties hold for any two distinct $i,j<n_{rp}$:
\begin{itemize}
    \item $O_i$ is an open subset of $[\mathbb{R}]^{n+1}$ and $W_i$ is an open subset of $[\mathbb{R}]^{|a_i|}$ (whenever $i\in P)$.
    \item $x_{n_r+i}\in O_i$ and $a_i\in W_i$ (for $i\in P)$,
    \item $\overline{W_i}<\overline{W_j}$ whenever $i<j\in P$,
    \item $\bigcup \overline{O_i}\cap \bigcup \overline{O_j}=\emptyset$ (So in particular, $\overline{O_i}\cap \overline{O_j}=\emptyset$),
    \item $\type(O_i,O_j)=\type(x_i,x_j),$
\end{itemize}
 We proceed to define by reverse recursion two more sequences $\langle U_i\rangle_{i<n}$ and $\langle V_i\rangle_{i\in P}$ inside $M$ satisfying the following conditions for $i<n$:
 \begin{enumerate}[label=(\Roman*)]
     \item  $U_i\subseteq O_i$, it is open in $[\mathbb{R}]^{n+1}$ and $x_{n_r+i}\cap \bigcup U_i=\emptyset,$
     \item $V_i\subseteq W_i$ and it is open in $[\mathbb{R}]^{|a_i|}$ and $a\cap \bigcup V_i=\emptyset$ whenever $i\in P$,
    \item $\type(x_{n_r+i},U_i)$ is not constant,
     \item $\Vert\type(a,\otimes_{i\in P} V_i)\Vert=\Vert \tau\Vert$,
     \item there are $q\in \mathbb{P}(\X,\E)\cap M$ and $b\in (\otimes_{j\in P}V_j)\cap M$ such that:\begin{itemize}
         \item $q$ end-extends $r$,
         \item $n_q=n_p$,
         \item $x^q_{n_r+i}\in U_i$ for each $i<n_{rp}$.
         \item $q\forces{$b\in \dot{A}$ }$.
     \end{itemize}
 \end{enumerate}
 It should be clear that for any $q$ and $b$ as in (V), the conditions (I) and (III) assure that $q$ is compatible with $p$. Furthermore, (IV) implies that $\Vert\type(a,b)\Vert=\Vert \tau \rVert$. Thus, the proof will be over once we finalize the construction of afore mentioned sequences.\\
 Before we start the construction, we name one last set of objects. Let $m_{n_{rp}+1}=|a|$ and given $i\leq n_{rp}$, let $m_i\leq |a|$ be such that $a[m_i]=a\cap \mathcal{M}_{n_{rp+i}} $. Also let $\sigma_i=: |a_i| \longrightarrow 2$ be the type defined as:$$\sigma_i(j)=\tau(m_i+j).$$
For each $i\leq n_{rp}$ let $p_i=r\cup\{(\mathcal{M}_{n_r+j},x_{n_r+j})\,:j<i\}$ and note that $p_i$ is a condition in $\mathbb{P}(\X,\E)$ which is end-extended by $p$ and such that $p_i\in \mathcal{M}_{i+1}(0).$ Furthermore, $p_{n_{rp}}=p$ and $p_0=r.$ Finally, let $\mathcal{D}$ be the set of all pairs $$(\langle y_i\rangle_{i<{n_{rp}}}, b)\subseteq(\prod_{i<n_{rp}} O_i)\times(\bigotimes_{i\in P}W)$$ for which there is $q\in \mathbb{P}(\X,\E)$ end-extending $r$ so that $n_q=n_p$, $x^q_{n_r+i}=y_i$ for each $i$, and $q\forces{ $a\in \dot{A}$ }$. Observe that $D\in \hat{M}$ and $D\in H(\theta)$. Thus, $D\in M=H(\theta)\cap \hat{M}$. Even more, $(\langle x_{{n_{r}+i}}\rangle_{i<{n_rp}},a)\in\mathcal{D}$.\\\\
\noindent
 The construction begins: Assume $k\leq n$ an we have defined  $\langle U_i\rangle_{k<i<n}$ and $\langle V_i\rangle_{k<i\in P}$ inside $M$ satisfying the properties (I), (II) and (III). Furthermore, assume  that following recursive hypotheses hold:\begin{enumerate}
     \item[(IV)$_{r}$] For any $k<i\in P$, $\lVert\type(a_i, V_i)\lVert=\lVert \tau\rVert$ if $a=a_i$, and $\type(a_i,V_i)=\sigma_i$ if $a\not=a_i.$
     \item[(V)$_{r}$] There are $u\in \prod_{k<i<n_{rp}} U_i$ and $v\in \bigotimes_{k<i\in P}V_i$ so that $(\langle x_i\rangle_{i\leq n_r+k}\frown u,a[m_{k+1}]\cup v)\in \mathcal{D}.$
 \end{enumerate}
If $k\in P$ and $a_k\not=a$, we define $\mathcal{D}_k$ as the set of all $(y,z)\in O_k\times W_k$ for which there are $u\in \prod_{k<i<n_{rp}} U_i$ and $v\in \bigotimes_{k<i\in P}V_i$ so that $$(\langle x_i\rangle_{i<n_r+k}\frown y\frown u,a[m_{k}]\cup z\cup v)\in \mathcal{D}.$$
It should be clear that $\mathcal{D}\in \mathcal{M}_{n_r+k}(0)$. Furthermore, $(x_{n_r+k},a_k)\in \mathcal{D}$. Thus, we can use Lemma \ref{lemmakind2} to find $(y,b)\in \mathcal{D}_k\cap \mathcal{M}_{n_r+k}(0)$ so that $\type(x_{n_r+k},y)$ is not constant and  $\type(a_k,z)= \sigma_k$. Now, we just enlarge $y$ and $z$ to open sets $U_k\subseteq O_k$ and $V_k\subseteq W_k$ inside $M$ whose union is disjoint from $x_{n_r+k}$ and $a_k$ respectively, and such that $\type(x,y)=\type(x,U_k)$ and $\type(a,z)=\type(a, V_k)$. It is straightforward that the recursive hypotheses still hold for $k$. Now, if $k\in P$ and $a_k=a$ we argue in a similar way but using Lemma \ref{lemmakind1}, and if $k\not\in P$ we do the same but by making use of Lemma \ref{lemmakind2empty}.
This finishes the construction.\\
We have constructed two sequences $\langle U_i\rangle_{i\leq n}$ and $\langle V_i \rangle_{i\in P}$ which clearly satisfy (I), (II), (III) and (V). The key to prove that  that they also satisfy (IV) is to note one of the two following cases occur: \begin{enumerate}
    \item  $P=\{i\}$ for some $i\leq n$, situation in which $a_i=a$ and consequently $\sigma_i=\tau$.  Thus $\Vert\type( a,V_i)\Vert=\Vert\type(a,\otimes_{i\in P} V_i)\Vert=\Vert\tau\Vert$.
    \item $|P|\geq 2$, which means that $a_i\not=a$ for each $i\in P$. Thus, $\type(a_i,V_i)=\sigma_i$ for each $i\in P$. As $a_{P(0)}<\dots <a_{P(|P|-1)}$ and $V_{P(0)}<\dots <V_{P(|P|-1)}$, we have that $\type(a, \otimes_{i\in P} V_i)=\sigma_{P(0)}\frown\dots\frown\sigma_{P(|P|-1)}=\tau.$
    \end{enumerate}
\end{proof}
\end{theorem}
As a corollary of the previous theorem, we get the main theorem of the paper.
\begin{theorem}[Previously enumerated as Theorem \ref{maintheorem}] Let $n\geq 2$. The statement \say{There is an $n$-entangled set, but there are no $n+1$-entangled sets}+PID+MA is consistent.
\end{theorem}

\section{Some preliminaries on construction schemes}\label{preliminariesschemesection}
We dedicate this section present some preliminaries on construction schemes. All material presented here is spread throughout \cite{finitizationclubch}, \cite{schemescruz}, \cite{forcingschemescruz} and \cite{schemenonseparablestructures}. The reader is referred to the author's Phd thesis \cite{chapital2024constructionschemesbuildinguncountable} for a complete introduction to this topic.\\
A \emph{(scheme) type} is a sequence $\tau=\langle m_k,n_{k+1},r_{k+1}\rangle_{k\in\omega}$ of triplets of natural numbers such that the following conditions hold for any $k\in\omega$:\\
\begin{enumerate}[label=$(\alph*)$,itemsep=0.5em]
\begin{minipage}{5cm}
\item $m_0=1,$
\item $n_{k+1}\geq 2,$
\end{minipage}
\begin{minipage}{5cm}
\item $m_k>r_{k+1}, $
\item $m_{k+1}=r_{k+1}+(m_k-r_{k+1})n_{k+1}.$
\end{minipage}\\
\end{enumerate}
A \emph{construction scheme\footnote{Sometimes, we will call it simply a scheme.} (of type $\tau$)} is a family $\mathcal{F}\subseteq \text{Fin}(\omega_1)$ which: Is cofinal in $(\text{Fin}(\omega_1),\subseteq )$, any member of $\mathcal{F}$ has cardinality $m_k$ for some $k\in\omega$, and furthermore, if we put $\mathcal{F}_k:=\{ F\in \mathcal{F}\,:\,|F|=m_k\}$, then the following two properties are satisfied for each $k\in\omega$:
\begin{enumerate}[label=(\roman*),itemsep=0.5em]
\item $\forall F,E\in \mathcal{F}_k\big(\; E\cap F\sqsubseteq E,F\;\big)$,
\item $\forall F\in \mathcal{F}_{k+1}\;\exists F_0,\dots,F_{n_{k+1}-1}\in \mathcal{F}_k$ such that $$F=\bigcup\limits_{i<n_{k+1}}F_i.$$
Moreover, $\langle F_i\rangle_{i<n_{k+1}}$ forms a root-tail-tail $\Delta$-system with root $R(F)$ such that  $|R(F)|=r_{k+1}$ and $R(F)<F_0\backslash R(F)<\dots < F_{n_{k+1}-1}\backslash R(F).$
\end{enumerate}
Note that for a given $F\in \mathcal{F}_{k+1}$, each of the $F_i$'s mentioned above can be written as $F[r_{k+1}]\cup F[\,[a_i,a_i+1)\,]$ where $a_i=r_{k+1}+i\cdot(m_{k+1}-m_k)$. In particular, this is saying that the family $\langle F_i\rangle_{i<n_{k+1}}$ is unique, so  we call it the \emph{canonical decomposition} of $F$. 

\begin{definition}[The $\rho$-function] Let $\mathcal{F}$ be a scheme. We define $\rho:\omega_1^2\longrightarrow \omega$ as:$$\rho(\alpha,\beta)=\min (k\in\omega\,:\,\exists F\in \mathcal{F}\,(\{\alpha,\beta\}\subseteq F )\,).$$
For each finite $A\subseteq \omega_1$, we also define $$\rho^A=\max(\rho[A^2])=\max(\rho(\alpha,\beta)\,:\,\alpha,\beta\in A).$$    
\end{definition}
It is not hard to see that $\rho^F=n$ for each $F\in \mathcal{F}_n$.
The most important feature of the function $\rho$ is that it is an ordinal metric. This means that it satisfies the properties stated in the following lemma.
\begin{lemma}Let $\mathcal{F}$ be a construction scheme. The following properties hold for any $\alpha,\beta,\gamma\in dom(\mathcal{F})$ and each $k\in\omega$:\vspace{0.5em}
\begin{enumerate}[label=$(om_{\arabic*})$, itemsep=0.5em]
\item $\rho(\alpha,\beta)=0$ if and only if $\alpha=\beta$.
\item $\rho(\alpha,\beta)=\rho(\beta,\alpha).$ 
\item If $\alpha\leq \min(\beta,\gamma)$, then $\rho(\alpha,\beta)\leq \max(\,\rho(\alpha,\gamma),\rho(\beta,\gamma)\,)$.
\item $\{\xi\leq\alpha\,:\,\rho(\alpha,\xi)\leq k\}$ is finite.
\end{enumerate}
\end{lemma}
Given $\alpha\in dom(\mathcal{F})$ and $k\in\omega$ we can define the $k$-closure of $\alpha$ as $(\alpha)_k:=\{\xi\leq \alpha\,:\rho(\alpha,\xi)\leq k\}$ and $(\alpha)_k^-:=(\alpha)_k\backslash \{\alpha\}$. Note that property $(om_4)$ is saying that all the $k$-closures are finite. It is a useful fact that for any $k\in\omega$ (with $m_k\leq |dom(\mathcal{F})|$) there is at least one $F\in \mathcal{F}_k$ such that $\alpha\in F$. Even more, for any such $F$ we have the equalities: $$F\cap (\alpha+1)=(\alpha)_k$$
$$F\cap \alpha=(\alpha)^-_k.$$
Given $\alpha\in dom(\mathcal{F})$, the $k$-cardinality function $\lVert \alpha\rVert_{\_}:\omega\longrightarrow \omega$ is defined as:
$$\lVert \alpha\rVert_k=|(\alpha)^-_k|.$$

 \begin{definition}[The $\Delta$-function]Let $\mathcal{F}$ be a scheme. We define $\Delta:dom(\mathcal{F})^2\longrightarrow \omega+1$ which is defined as: $$\Delta(\alpha,\beta)=\Delta(\lVert \alpha\rVert_{\_},\lVert \beta\rVert_{\_}).$$
 \end{definition}
\begin{definition}[The $\Xi$-function]Let $\mathcal{F}$ be a scheme. Given   
$\alpha\in \omega_1$ we define the function $\Xi_\alpha:\omega\longrightarrow \{-1\}\cup \omega$ as follows;  If $1\leq k\in\omega$, then there is $F\in \mathcal{F}_k$ so that $\alpha\in F$. According to the property $(ii)$ in the definition of a construction scheme, we have that either $\alpha\in R(F)$ or there is a unique $i<n_k$ such that $\alpha\in F_i\backslash R(F)$. We then define $$\Xi_\alpha(k):=\begin{cases}-1&\text{ if }\alpha\in R(F)\\
i&\text{ if }\alpha\in F_i\backslash R(F)\end{cases}$$
It can be proved that this definition does not depend on the choice of $F.$ Now, if $k=0$, we simply define $\Xi_\alpha(k)$ as $0$.
\end{definition}

The next lemma relates the functions $\rho$, $\Xi$ and $\Delta$ in a very nice way. 
\begin{lemma}\label{lemmaxi}Let $\mathcal{F}$ be a construction scheme, $\alpha<\beta\in dom(\mathcal{F})$ and $1\leq k\in \omega$ be such that $m_k\leq |dom(\mathcal{F})|$. Then:\vspace{0.5em}
\begin{enumerate}[label=$(xi_{\alph*})$, itemsep=0.5em]
\item If $k<\Delta(\alpha,\beta)$, then  $\Xi_\alpha(k)=\Xi_\beta(k).$
\item If $k=\rho(\alpha,\beta)$, then $0\leq \Xi_\alpha(k)<\Xi_\beta(k).$
\item If $k>\rho(\alpha,\beta)$, then either $\Xi_\alpha(k)=-1$ or $\Xi_\alpha(k)=\Xi_\beta(k).$ 
\item If $k=\Delta(\alpha,\beta)$ then $0\leq \Xi_\alpha(k)\not=\Xi_\beta(k)\geq 0.$ 
\end{enumerate}
\end{lemma}
Note that another way of stating $(xi_a)$ and $(xi_d)$ is simply by saying that $\Delta(\alpha,\beta)=\Delta(\Xi_\alpha,\Xi_\beta)$.
\begin{definition}[Captured families]Let $\mathcal{F}$ be a scheme, $1\leq n,m\in \omega$ and $\mathcal{D}=\langle D_i\rangle_{i<n}\subseteq[\omega_1]^m$. Given $l\in \omega$, we say that $\mathcal{D}$ is $\Delta$-captured at level $l$ if:
\begin{itemize}
    \item $\mathcal{D}$ forms a root-tail-tail $\Delta$-system with some root $R$ of cardinality $r$.
    \item For all  $a<m$  and $j<n$:
    $\Xi_{D_j(a)}(l)=\begin{cases}-1&\text{ if }i<r\\
    \:j&\text{ if }i\geq r
    \end{cases}$
    \item $\Delta(D_j(a),D_i(a))=l$ for any two distinct $i,j<n$ and every $r\leq a<m$.
\end{itemize}
If moreover $\rho(D_j(a),D_i(a))=l$ for any two distinct $i,j<n$ and every $r\leq a<m$, we say that $\mathcal{D}$ is \emph{captured.} Given a finite $D\subseteq \omega_1$, we say that $D$ is $\Delta$-captured (captured) whenever $\langle \{D(i)\}\rangle_{i<|D|}$ is $\Delta$-captured (captured).
\end{definition}
\begin{definition} Let construction scheme $\mathcal{F}$ is said to be:
\begin{itemize}
    \item $2$-$\Delta$-capturing (resp. $2$-capturing) if for each uncountable $\mathcal{S}\subseteq \Fin(\omega_1) $ there are infinitely many $l\in \omega$ for which there is $\mathcal{D}\in [\mathcal{S}]^2$ which is $\Delta$-captured (res. captured) at level $l$.
    \item \emph{Fully $\Delta$-capturing} (resp. fully capturing) if for each uncountable $\mathcal{S}\subseteq \Fin(\omega_1) $ there are infinitely many $l\in \omega$ for which there is $\mathcal{D}\in [\mathcal{S}]^{n_l}$ which is $\Delta$-captured (res. captured) at level $l$.
\end{itemize}
\end{definition}
\begin{proposition}\label{capturingequiv}A scheme $\mathcal{F}$ is $2$-$\Delta$-capturing if for each uncountable $\mathcal{S}\subseteq [\omega_1]^2$ family of pairwise disjoint sets there are infinitely many $l\in \omega$ for which there is $\mathcal{D}\in [\mathcal{S}]^2$ which is $\Delta$-captured  at level $l$. 
\end{proposition}
\begin{proposition}\label{deltacapturingequiv}A scheme $\mathcal{F}$ is $2$-capturing if for each uncountable $S\subseteq \omega_1$there are infinitely many $l\in \omega$ for which there is $\mathcal{D}\in [S]^2$ which is $\Delta$-captured  at level $l$. 
\end{proposition}
The Full capturing axiom (FCA) is the statement \say{There is a fully capturing scheme for any type}.  We define in an analogous way the Fully $\Delta$-capturing axiom (FCA$^\Delta$). In \cite{forcingandconstructionschemes}, it was proved that FCA holds in the Cohen model. We proved the following theorem in \cite{schemescruz}. 
\begin{theorem}[Under $\Diamond$,\cite{schemescruz}]FCA holds
\end{theorem}
In \cite{finitizationclubch}, we proved the analogous result for fully-$\Delta$-capturing schemes but assuming only CH.
\begin{theorem}[Under CH,\cite{finitizationclubch}]$FCA^\Delta$ holds.
\end{theorem}
For a scheme $\mathcal{F}$, e say that a forcing notion $\mathbb{P}$ \emph{$2$-preserves (resp $2$-$\Delta$-preserves)} $\mathcal{F}$ whenever $\mathbb{P}$ forces $\mathcal{F}$ to be $2$-capturing (resp. $2$-$\Delta$-capturing).  We define $\mathfrak{m}_\mathcal{F}$ as the minimum of all the Martin's numbers associated to $ccc$-forcing notions which $2$-prerserve $\mathcal{F}$. Analogously $\mathfrak{m}^\Delta_\mathcal{F}$ denotes the minimum of all the Martin's numbers associated to $ccc$-forcing notions which $2$-$\Delta$-prerserve $\mathcal{F}$. MA$(\mathcal{F})$(resp. MA$^\Delta(\mathcal{F})$) denote the statement $\mathfrak{m}_\mathcal{F}=\mathfrak{c}>\omega_1$ (resp.$\mathfrak{m}^2_\mathcal{F}=\mathfrak{c}>\omega_1$).

In \cite{forcingschemescruz}, we introduced and studied thoroughly the cardinal $\mathfrak{m}_\mathcal{F}$. There we proved:
\begin{theorem}(\cite{forcingschemescruz}) Given a $2$-capturing scheme $\mathcal{F}$, there is an extension of the universe via a $ccc$-forcing in which MA$(\mathcal{F})$ holds.
\end{theorem}

It is easy to see that the exact same arguments yield an analogous theorem for $2$-$\Delta$-capturing schemes. The main tools for showing that a forcing $\mathbb{P}$ is $ccc$ and $2$-preserves (resp. $2$-$\Delta$-preserve a scheme) are implicit in the equivalences presented above in propositions \ref{capturingequiv} and \ref{deltacapturingequiv}.
The following lemma is proved in \cite{forcingschemescruz}.

\begin{lemma}\label{lemmaequivalencefunctioncapturingpreserving}Let $\mathcal{F}$ be an $2$-capturing  scheme and $\mathbb{P}$ be a forcing. The two following statements are equivalent:
\begin{enumerate}[label=$(\arabic*)$]
    \item $\mathbb{P}$ is $ccc$ and $n$-preserves $\mathcal{F}.$
    \item For any $\mathcal{A}\in[\mathbb{P}]^{\omega_1}$ and each injective function $\nu:\mathcal{A}\longrightarrow \omega_1$, there are distinct $p,q\in \mathcal{A}$ for which:
    \begin{itemize}
        \item $p$ and $q$ are compatible. 
        \item $\{\eta(p),\eta(q)\}$ is captured.
    \end{itemize}
\end{enumerate}
\end{lemma}
In a complete analogous way we can prove the following lemma. We omit the proof.
\begin{lemma}\label{lemmaequivalencefunctiondeltacapturingpreserving}Let $\mathcal{F}$ be an $2$-$\Delta$-capturing  scheme and $\mathbb{P}$ be a forcing. The two following statements are equivalent:
\begin{enumerate}[label=$(\arabic*)$]
    \item $\mathbb{P}$ is $ccc$ and $2$-$\Delta$-preserves $\mathcal{F}.$
    \item For any $\mathcal{A}\in[\mathbb{P}]^{\omega_1}$ and each function $\nu:\mathcal{A}\longrightarrow [\omega_1]^2$ with pairwise disjoint images, there are distinct $p,q\in \mathcal{A}$ for which:
    \begin{itemize}
        \item $p$ and $q$ are compatible. 
        \item $\{\eta(p),\eta(q)\}$ is $\Delta$-captured.
    \end{itemize}
\end{enumerate}
\end{lemma}
\section{$n$-entangled sets and full $\Delta$-capturing}\label{schemenotn+1section}
The main goal of this section is to show that the existence of a fully $\Delta$-capturing scheme (with a a fastly growing type) yields the existence of $n$-entangled sets which are not $n+1$-entangled. Rather than constructing such entangled sets over the reals, we will do it over $\omega^\omega$ endowed with the order $f<_{lex} g$ if $f(n)<g(n)$ where $n=\Delta(f,g)$. It is easily seen that this total order is separable, hence it can be embedded in $\mathbb{R}.$ Now, the entangled sets which we construct will be indexed by countable ordinals. For that reason, we will need the equivalence presented below which the author assume is folkore (see \cite{finitizationclubch} for a proof).
\begin{lemma} Let $(\mathbb{X},<)$ be a total order, $k\in\omega$ and $\E\in [\mathbb{X}]^{\omega_1}$ injectively enumerated as $\langle e_\alpha\rangle_{\alpha\in \omega_1}$. Then $\mathcal{E}$ is $n$-entangled if and only if for every uncountable family $\mathcal{A}\subseteq [\omega_1]^{n}$ of pairwise disjoint sets and  each $t:n\longrightarrow 2$ there are distinct $c,d\in \mathcal{C}$ for which $\Vert\bar\type(C,D)\Vert=\Vert\tau\Vert$. Here, $\bar\type(C,D):n\longrightarrow 2$ is defined by $\bar\type(C,D)(i)=0$ if and only if $e_{C(i)}<e_{D(i)}.$
\end{lemma}

For the rest of this section, we fix $n\in \omega$ and a fully $\Delta$-capturing scheme $\mathcal{F}$ with type $\langle m_k,n_{k+1},r_{k+1}\rangle_{k\in \omega}$ in such way that $n_{k+1}\geq 2^{m_k}$ for each $k\in \omega$. For any such $k$, let us enumerate (possibly with repetitions) the set $[m_k\backslash r_{k+1}]^{n}$ as $\langle  C^k_i\rangle_{0< i<n_{k+1}}$.  Given $\alpha\in \omega_1$, we define two functions:\begin{enumerate}
    \item $e^n_\alpha:\omega\longrightarrow \mathbb{Z}$ given by:
   $$e^n_\alpha(k)=\begin{cases}0&\text{ if }\Xi_\alpha(k)\leq 0\\
   \Xi_\alpha(k)&\text{ if }\Xi_\alpha(k)>0\text{ and }\lVert \alpha\rVert_{k-1}\in C^{k-1}_{\Xi_\alpha(k)}\\
   \Xi_\alpha(k)&\text{ if }\Xi_\alpha(k)>0,\,|C^{k-1}_{\Xi_\alpha(k)}|=n\text{ and }\Vert \alpha\Vert_{k-1}<\min(C^{k-1}_{\Xi_\alpha(k)})\\
   
  - \Xi_\alpha(k)&\text{ in any other case }
  \end{cases}$$
  
    \item $o^n_\alpha:\omega\longrightarrow \mathbb{Z}$ given by:
   $$o^n_\alpha(k)=\begin{cases}0&\text{ if }\Xi_\alpha(k)\leq 0\\
   \Xi_\alpha(k)&\text{ if }\Xi_\alpha(k)>0\text{ and }\lVert \alpha\rVert_{k-1}\in C^{k-1}_{\Xi_\alpha(k)}\\
  - \Xi_\alpha(k)&\text{ if }\Xi_\alpha(k)>0\text{ and }\Vert \alpha\Vert_{k-1}\notin C^{k-1}_{\Xi_\alpha(k)}
   \end{cases}$$
\end{enumerate}
Since $n$ is fixed, let us write for now $e^n_\alpha$ and $o^n_\alpha$ simply as $e_\alpha$ and $o_\alpha$ respectively. Also let $\mathbb{E}_e=\mathbb{E}^n_e:=\{e_\alpha\,:\,\alpha\in \omega_1\}$ and $\mathbb{E}^n_e=\mathbb{E}^n_o:=\{o_\alpha\,:\,\alpha\in \omega_1\}$.  We proceed to analyze the degree of entangledness of these two sets.
\begin{proposition}\label{entangledeven}$\mathbb{E}_e$ is $2n$ entangled.
    \begin{proof}Let $\tau:2n\longrightarrow 2$ be a type and $\mathcal{A}\subseteq [\omega_1]^{2n}$ be an uncountable family with pairwise disjoint elements. By changing $\tau$ by $\tau^1$ if necesarry, we may assume that $\tau^{-1}[\{0\}]\leq n$. Furthermore, if $\tau^{-1}[\{0\}]=n$ we may also assume that $\tau(0)=0$. As $\mathcal{F}$ is fully $\Delta$-capturing, there are $l\in \omega$ and a family $\{A_0,\dots,A_{n_l-1}\}\in [\mathcal{A}]^{n_l}$ which is $\Delta$-captured at level $l$. Let $$C:=\{\lVert A_0(i)
    \Vert_{l-1}\,:\,i<2n\text{ and }\tau(i)=0\}$$
    Then $C\subseteq m_{l-1}\backslash r_l$ and $|C|\leq n$. In this way, there is $j<n_l$ so that $C=C^{l-1}_j$. We claim that $\Vert \bar\type(D_0,D_j)\Vert=\Vert\tau\Vert$. Indeed, let $i<2n$.  Then:
    \begin{enumerate}
        \item $e_{D_0(i)}(l)=0$ because $\Xi_{D_0(i)}(l)=0$,
        \item $\lVert D_0(i)\Vert_m=\Vert D_j(i)\Vert_m$ for any $m<l$ because $\Delta(D_0(i),D_j(i))=l$.
    \end{enumerate} From the two points above it easily follows that $e_{D_0(i)}(m)=e_{D_j(i)}(m)$  for any such $m$.  In order to finish, just observe that:
    \begin{itemize}
        \item If $\tau(i)=0$, then $j=\lVert D_j(i)\Vert_{l-1}=\Vert D_0(i)\Vert_{l-1}\in C=C^{l-1}_j$. Thus, $e_{D_j(i)}(l)=j>0=e_{D_0(i)}(l)$.
        \item If $\tau(i)=1$, then $j=\lVert D_j(i)\Vert_{l-1}=\Vert D_0(i)\Vert_{l-1}\notin C=C^{l-1}_j$. In this way, $e_{D_j(i)}=-j<0=e_{D_0(i)}(l)$ (If $|C|<n$, this is clear by definition. If $|C|=n$, this is because we assume that $\tau(0)=0$).
    \end{itemize}    
    \end{proof}
\end{proposition}
In a completely similar way, we can prove the following result.
\begin{proposition}$\mathbb{E}_o$ is $2n+1$ entangled.
\end{proposition}
We proceed to bound from above the maximum amount of entangledness of both $\mathbb{E}_e$ and $\mathbb{E}_o$. What we will do is to take an uncountable family of pairwise disjoint intervals over $\omega_1$ and show that we can refine that family in order to avoid the alternating types (those for which  any two numbers have the same value).

\begin{proposition}$\mathbb{E}_e$ is not $2n+1$ entangled.
\begin{proof}For each $\gamma\in \Lim$, let $D_\gamma:=[\gamma,\gamma+2n+1)=(\gamma+2n+1)\backslash \gamma$. By the pigeoholeprinciple, we can find $X\in [\omega_1]^{\omega_1}$ and $k,a\in \omega$ so that $\rho^{D_\gamma}=k$ for each $\gamma\in X$ and $\lVert D_\gamma(0)\Vert_k=\Vert\gamma\Vert_k=a$. We claim that for any two distinct $\delta,\gamma\in X$ we have that $\Vert\bar \type(D_\delta,D_\gamma)\Vert\not=\tau$ where $\tau:2n+1\longrightarrow 2$ is alternating type defined by: $$\tau(i)=\begin{cases}
    0&\text{ if i is even}\\
    1&\text{ if i is odd}
\end{cases}$$
Assume towards a contradiction that this is not true. That is, there are $\gamma,\delta\in X$ with $\Vert\bar \type(D_\delta,D_\gamma)\Vert=\tau$. Let $l=\Delta(\gamma,\delta)$ and note that since $\Vert \gamma\Vert_k=\Vert\delta\Vert_k$, then $l>k$.\\\\
\underline{Claim 1}: For any $i<2n+1$, $\Delta(\gamma+i,\delta+i)=l$.
\begin{claimproof}[Proof of claim]Note that for any $k<s\in \omega$, then there is $F\in \mathcal{S}$ with $[\gamma,\gamma+2n+1)\subseteq F$. By definition, we have that $$\lVert \gamma+i\Vert_s=|(\gamma+i)^-_s|=|F\cap( \gamma+i)|=|F\cap \gamma|+i=\Vert\gamma\Vert_s+i.$$
and in a similar way we get that $\lVert \delta\lVert_s+i=\lVert\delta+i\lVert_s$. The claim follows directly from these inequalities.  
\end{claimproof}
From the arguments presented in  Proposition \ref{entangledeven}, we conclude that $\Delta(e_{\gamma+i},e_{\delta+i})=\Delta(\gamma+i,\delta+i)=l$ for any $i<2n+1$. This means that the order between  $e_{\gamma+i}$ and $e_{\gamma+i}$ is determined by the value of these two functions at $l$.\\
Observe that  $\rho(\gamma+i,\gamma)\leq \rho^{D_\gamma}=k<l$ and $\rho(\delta+i,\delta)\leq \rho^{D_\delta}=k<l$ for each $i<2n+1$. In this way, $$0\leq \Xi_\gamma(l)=\Xi_{\gamma+j}(l)\not=\Xi_{\delta+j}(l)=\Xi_\delta(l)\geq 0.$$
Without any loss of generality let us assume that $j:=\Xi_\gamma(l)>\Xi_\delta(l)=:z$. Observe that for each $i$ we have that $e_{\gamma+i}(l)<e_{\delta+i}(l)$ if and only if $e_{\gamma+i}(l)=-j$. Let $C_{even}=\{\lVert \gamma+i\Vert_{l-1}\,:\,i<2n+1\text{ is even}\}$ and $C_{odd}=\{\{\lVert \gamma+i\Vert_{l-1}\,:\,i<2n+1\text{ is odd}\}$. Consider the following cases:\\\\
\underline{Case 1}:  If $e_\gamma(l)=j$ then $e_{\gamma+i}(l)=j$ for any $i\in C_{even}$. Since we are assuming that the type is alternating, then $e_{\gamma+1}(l)=-j$. This implies that $\lVert \gamma+1\Vert_{l-1}>\min(C^{k-1}_j)$. Consequently, $\Vert \gamma+1\Vert-1=\lVert \gamma\Vert_{l-1}\geq\min(C^{k-1}_j)$. This means that $C_{even}\subseteq C^{l-1}_j$. However, this is a contradiction  since $|C_{even}|=n+1>n\geq |C^{l-1}_j|$.\\\\
\noindent
\underline{Case 2}: If $e_\gamma(l)=-j$,  we can argue in a similar way to conclude that $C_{odd}\subseteq C^{l-1}_j$. If $|C^{l-1}_j|<n$, this leads to the same contradiction as in the previous case. On the other hand, if $|C^{l-1}_j|=n$, then the fact that $e_\gamma(l)=-j$ implies that $\min(C^{l-1}_j)<C_{odd}$. In this way,  $n=|C_{odd}|<|C^{l-1}_j|-1=n-1$ which is again a contradiction.
\end{proof} 
\end{proposition}
By arguing in a similar way but ignoring Case 1, we get:
\begin{proposition} $\mathbb{E}_o$ is not $2n+2$ entangled.
\end{proposition}
We can use the same ideas as above to prove that the family $\mathbb{E}_{\mathcal{F}}=\langle \Xi_\alpha\rangle_{\alpha\in \omega_1}$ is a subset of $\omega^\omega$ which is not even $2$-entangled.
Through out the proofs of this section we have been   using the fact that $\Delta(o_\alpha,o_\beta)=\Delta(\alpha,\beta)=\Delta(e_\alpha,e_\beta)$ for any two $\alpha,\beta\in \omega_1$. Since $\Delta(\alpha,\beta)=\Delta(\Xi_\alpha,\Xi_\beta)$ as well (see Lemma \ref{lemmaxi}), this means that $\mathbb{E}^\mathcal{F}$ is homeomorphic to both $\mathbb{E}_e$ and $\mathbb{E}_o$ via the functions sending $\Xi_\alpha$ to $e_\alpha$ and $o_\alpha$ respectively. In \cite{finitizationclubch}, we defined a set $\mathbb{E}_\infty$ which is entangled and homeomorphic to $\mathbb{E}^\mathcal{F}$ for the same reasons. Thus, the results of this section together to the one from \cite{finitizationclubch} yield following theorem generalizing the one from \cite{carroy2025questionsentangledlinearorders}.

\begin{theorem}[Under $FCA^\Delta$(Previously enumerated as Theorem \ref{schemenotn+1theorem}] There is a family $\langle \mathbb{E}_n\rangle_{n\in \omega}$ of pairwise homeomorphic uncountable subsets of $\mathbb{R}$ so that $\mathbb{E}_0$ is entangled, $\mathbb{E}_1$ is not $2$-entangled and $\mathbb{E}_n$ is an $n$-entangled but not $n+1$-entangled set for any $n\geq 2.$
    
\end{theorem}

\section{More on schemes and entangled sets}\label{mfsection}
 As we have mentioned before, the existence of fully $\Delta$-capturing scheme with a rapidly growing (scheme) type implies the existence of an entangled set. In particular, this means that the existence of certain fully capturing schemes imply the same result. Thus, it is natural to ask what is the amount of capturing needed in order to build this kind of objects. The main goal of this section is to show that neither the existence of a $2$-$\Delta$-capturing schemes or even $2$-capturing schemes implies the existence of an entangled set. In order to understand the importance of these results it is worth to remark that most of the objects which we know how to construct using schemes just require the scheme is $2$-capturing (see \cite{schemescruz}). For example, in \cite{schemescruz} we showed that the existence of a $2$-capturing scheme implies the failure of BA$(\omega_1)$. Actually, it is easy to prove that $\{\Xi_\alpha\:\,\alpha\in \omega_1\}$ is an increasing set\footnote{A total order $\mathbb{X}$ is said to be increasing if for any $1\leq n\in \omega$ and each uncountable pairwise disjoint family $\mathcal{A}\subseteq[\X]^n$, there are distinct $A<B\in \mathcal{A}$ so that $\type(A,B)$ is constant. The existence of increasing sets also implies the failure of BA$(\omega_1)$. These objects were studied for the first time in \cite{MAdoesnotImplyBA}.} by just assuming that $\mathcal{F}$ is $2$-$\Delta$-capturing. 
\begin{theorem}[Under $\mathfrak{m}_\mathcal{F}>\omega_1$(Previously enumerated as Theorem \ref{mFtheorem})]There are no entangled sets.
\begin{proof} Suppose towards a contradiction that there exists an entangled set $\mathbb{E}$. Fix an uncountable family of pairwise disjoint unordered pairs of $\mathbb{E}$ injectively enumerated as $\langle e_\alpha\rangle_{\alpha\in \omega_1}$. Let us define $\mathbb{P}$ to be the forcing consisting of all $A\in [\omega_1]^{<\omega}$ so that for any two $\alpha,\beta\in A$, $e_\alpha(0)<e_\beta(0)$ if and only if $e_\alpha(1)<e_\beta(1)$. In other words, $\type(e_\alpha,e_\beta)$ is constant. We order $\mathbb{P}$ by reverse inclusion. It is easy to see that $\mathbb{P}$ is a ccc-forcing provided that $\mathbb{P}$ is an entangled (thus, increasing) set. Furthermore, if $G\subseteq \mathbb{P}$ is an uncountable filter, then $\{e_\alpha\,:\alpha\in\bigcup G\}$ is an uncountable family of pairs of $\mathbb{E}$ such that $\type(e_\alpha,e_\beta)$ is constant for any two $\alpha,\beta\in \bigcup G$. The existence of any such family would contradict the fact that $\mathbb{E}$ is 2-entangled. Thus, since $\mathfrak{m}_\mathcal{F}>\omega_1$, it must be the case that $\mathbb{P}$ does not (2-)preserve $\mathcal{F}$. In this way, there exists an uncountable $\mathcal{A}\subseteq \mathbb{P}$ and an injective function $\zeta:\mathcal{A}\longrightarrow \omega_1$ so that for all $A,B\in \mathcal{A}$, if $\{\zeta(A),\zeta(B)\}$ is captured, then $A$ and $B$ are incompatible in $\mathbb{P}$. Without any loss of generality we may assume that all elements of $\mathcal{A}$ are pairwise disjoint and share a single cardinality, say $n$.  Furthermore, note that two elements of $\mathbb{P}$, say $A$ and $B$ are incompatible if and only if there are $\alpha\in A\backslash B$ and $\beta\in B\backslash A$ so that $\{\alpha,\beta\}\not\in\mathbb{P}$. In this way, we can also assume that the elements of $\mathcal{A}$ are pairwise disjoint.

Given $A\in \mathcal{A}$, define $E_A=\{e_\alpha(i)\,:\,\alpha\in A\text{ and }i\in 2\}$. Note that $A$ and $B$ are incompatible in $\mathbb{P}$, then $\type(E_A,E_B)$ is not constant.
We proceed to define  $\mathbb{Q}$ as the set of all $p\in [\mathcal{A}]^{<\omega}$  so that $\type(E_A,E_B)$ is not constant for any two distinct $A,B\in p.$\\\\
\noindent
\underline{Claim 1}: $\mathbb{Q}$ is $ccc$ and $2$-preserves $\mathcal{F}$.
\begin{claimproof}Let $\mathcal{B}\in[\mathbb{Q}]^{\omega_1}$ and $\nu:\mathcal{B}\longrightarrow \omega_1 $ be an injective function. We shall prove that there are two compatible $p,q\in \mathcal{B}$ so that $\{\nu(p),\nu(q)\}$ is captured. By refining $\mathcal{B}$, we may assume that there is $m\in \omega$ so that $\mathcal{B}\subseteq[\mathcal{A}]^n$, and by a similar argument as above, we may assume that the elements of $\mathcal{B}$ are pairwise disjoint. Let us enumerate each $p\in \mathcal{B}$ as $\{A_0^p,\dots A_{m-1}^p\}$. Given $i<m$, note that $E^p_i:=E_{A^p_i}\in[\mathbb{R}]^{2n}$.  Since being disjoint and having a non-constant defines an open graph, we can pick a sequence of disjoint open sets $\langle U_i^p\rangle_{i<m}\,:\,i\subseteq [\mathbb{R}]^{2n}$ inside a fixed countable base in such way that for any two distinct $i,j<m$:
\begin{itemize}
    \item $E^p_i\in U^p_i$,
    \item $\type(a,b)$ is not constant for any $a\in U^p_i$ and $b\in U^p_j$.
\end{itemize}
Since we can only choose from a countable collection of open sets, there is $\mathcal{C}\in [\mathcal{B}]^{\omega_1}$ and there are $U_0,\dots U_{m-1}$ so that $U^p_i=U_i$ for any $p\in \mathcal{C}$ and $i<m$. Given $p\in \mathcal{C}$, let $D_p=\{\zeta(A_i^p)\,:\,i<m\}\cup\{\nu(p)\}$. Again, there are $\mathcal{D}\in [\mathcal{C}]^{\omega_1}$, $d\in \omega$, $a\in d$ and $M\subseteq d$ so that for any  $p,q\in \mathcal{D}$ and every $i<m$:
\begin{enumerate}
\item $D_p\cap D_q=\emptyset,$
    \item $|D_p|=d$,
    \item $D_p(a)=\nu(p)$,
    \item $D_p(M(i))=\zeta(M(i)).$
\end{enumerate}
As $\mathcal{F}$ is $2$-capturing, it follows that there are distinct $p,q\in \mathcal{D}$ so that $\{D_p,D_q\}$ is captured. According to the item (1) above, we also have that $\{D_p(i),D_q(i)\}$ is captured for each $i<d$. In particular:
\begin{itemize}
\item$\{D_p(a),D_q(a)\}=\{\eta(p),\eta(q)\}$ is captured,
\item $\{D_p(M(i)),D_q(M(i))\}=\{\zeta(A_i^p),\zeta(A_i^q)\}$ is captured for any $i<m$.
\end{itemize}
By the first point above, the claim will be over once we prove that $p$ is compatible with $q$. Indeed, let $i,j<m$. If $i=j$, then the second point above implies that $\type(A^p_i,A^q_j)$ is not constant. On the other hand, if $i\not=j$, then $\type(A^p_i,A^q_j)$ is not constant because $A^p_i\in U_i$ and $A^q_j\in U_j$.
\end{claimproof} 
By the previous claim and since $\mathfrak{m}_\mathcal{F}>\omega_1$, there is an uncountable filter $G\subseteq \mathbb{Q}$. Note that there is no pair of elements inside  $\{E_A\,:\,\exists p\in G\,(A\in p)\}$ whose type is constant.This contradicts the fact that $\mathbb{E}$ is entangled.
\end{proof}
\end{theorem}

When dealing with $2$-$\Delta$-capturing schemes, we can actually prove a much more stronger result. This slighly improves the result from \cite{abrahamrubinshelah} stating that the existence of an increasing set is compatible with MA.
\begin{theorem}[Under $\mathfrak{m}^\Delta_\mathcal{F}>\omega_1$ (Previously enumerated as Theorem \ref{mfdeltatheorem})] Assume $\mathbb{P}$ is a $ccc$-forcing. Then $\mathbb{P}$ $2$-$\Delta$-preserves $\mathcal{F}$.
\begin{proof}Assume towards a contradiction that $\mathbb{P}$ does not $2$-$\Delta$-preserve $\mathcal{F}$. Therefore, there is $\mathcal{A}\in[\mathbb{P}]^{\omega_1}$ a function $\zeta:\mathcal{A}\longrightarrow [\omega_1]^2$ with pairwise disjoint images and with the property that for all  distinct $p,q\in \mathcal{A}$, if $\{\zeta(p),\zeta(q)\}$ is $\Delta$-captured, then $p$ and $q$ are incompatible. Let $\mathcal{X}=\zeta[\mathcal{A}]$.   We define $\mathbb{Q}$ as the forcing of all $p\in [\mathcal{X}]^{<\omega}$ so that $\{A,B\}$ is $\Delta$-captured for any two $A,B\in p.$ We order $\mathbb{Q}$ by reverse inclusion.\\\\
\noindent
\underline{Claim 1}: $\mathbb{Q}$ is $ccc$ and $2$-$\Delta$-preserves $\mathcal{F}$.
\begin{claimproof}Let $\mathcal{B}\in[\mathbb{Q}]^{\omega_1}$ and $\nu:\mathcal{B}\longrightarrow [\omega_1]^2$ be a function with pairwise disjoint images. By refining $\mathcal{B}$ we may assume that all of its elements are pairwise disjoint and have the same cardinality, say $m$. In this way, we enumerate each $p\in \mathcal{B}$ as $\{A^p_0,\dots,A^p_{n-1}\}$. We shall prove that there are $p,q\in \mathcal{B}$ so that $p$ and $q$ are compatible and $\{\nu(p),\nu(q)\}$ is $\Delta$-captured. Given $p\in \mathcal{B}$ let $D_p=\bigcup p\cup \nu(p)$ and $k_p=\rho^{D_p}$. Since $\mathcal{B}$ is uncountable, we can find $\mathcal{C}\in [\mathcal{B}]^{\omega_1}$, $k,d\in \omega$ and $M_0,\dots,M_{n-1},a\in[d]^{2}$ in such way that for any $p,q\in \mathcal{D}$ and every $i<m$:
\begin{enumerate}
    \item $D_p\cap D_q=\emptyset,$
    \item $|D_p|=d$ and $k_p=k$,
    \item $D_p[a]=\eta(p)$,
    \item $D_p[M_i]=A^p_i.$
\end{enumerate}
Since $\mathcal{F}$ is $2$-$\Delta$-capturing, can find distinct $p,q\in \mathcal{C}$ so that $\{D_p,D_q\}$ is $\Delta$-captured at some level $l>k$. In particular:
\begin{itemize}
    \item $\{D_p[a],D_q[a]\}=\{\nu(p),\nu(q)\}$ is $\Delta$-captured,
    \item $\{D_p[M_i],D_q[M_i]\}=\{A^p_i,A^q_i\}$ is $\Delta$-captured for any $i<m$.
\end{itemize}
In order to finish, it is enough to prove that $\{A^p_i,A^q_j\}$ is $\Delta$-captured for any two distinct $i,j<m.$ For this, first notice that, by definition, $\{A^p_i,A^p_j\}$ is $\Delta$-captured at some level $s$. Without any loss of generality, we may assume that
$\Xi_{(A^p_i(0))}(s)=\Xi_{A^p_i(1)}(s)=0$ and $\Xi_{(A^p_j(0))}(s)=\Xi_{A^p_j(1)}(s)=1$. By definition of $\Delta$-capturing we also have that $$\Delta(A^p_i(0),A^q_j(0))=s=\Delta(A^p_i(1),A^p_j(1)).$$
In particular, this means that $s=\Delta(A^p_i(0),A^p_j(0))\leq \rho(A^p_i(0),A^p_j(0))\leq \rho^{D_p}=k<l$. With this information, we are able to conclude that $\{A^p_i,A^q_j\}$ is also $\Delta$-captured at level $s$. Indeed, given $e\in 2$, we have that $\Delta(A^p_j(e),A^q_j(e))=l>s=\Delta(A^p_i(e),A^p_j(e))$. Thus:
\begin{itemize}
    \item $\Delta(A^p_i(0),A^q_j(0))=s=\Delta(A^p_i(1),A^q_j(1))$.
    \item $\Vert A^p_j(e)\Vert_s=\Vert A^q_j(e)\Vert_s$ for each $e\in 2$. This implies that, $\Xi_{A^q_j(0)}(0)=\Xi_{A^q_j(0)}(1)=1.$
\end{itemize}
\end{claimproof}
Since $\mathfrak{m}^\Delta_\mathcal{F}>\omega_1$, there is an uncountable $G\subseteq \mathbb{Q}$. Note that $\bigcup G$ is an uncountable subset of $\mathcal{X}$ such that $\{A,B\}$ is $\Delta$-captured for any two $A,B\in \bigcup G$. In this way,  we conclude that $\zeta^{-1}[\bigcup G]$ is an uncountable antichain in $\mathbb{P}$. This contradiction finishes the proof.
\end{proof}   
\begin{corollary}MA$^\Delta(\mathcal{F})$ implies MA.
\end{corollary}

\end{theorem} 
\bibliographystyle{plain}
\bibliography{bibliografia}
{\color{white}hlj}\\\\
\noindent
Jorge Antonio Cruz Chapital\\
Department of Mathematics, University of Toronto, Canada\\
cruz.chapital at utoronto.ca\\

\end{document}